\documentclass[leqno,12pt]{article}
\pdfoutput=1
\usepackage{latexsym}
\usepackage{amssymb,amsbsy,amsmath,amsfonts,amssymb,amscd}
\usepackage{graphicx}
\usepackage{color}
\usepackage{hyperref}
\usepackage{dsfont}

\newcommand {\e}  {\varepsilon}

\newtheorem{theorem}{Theorem}[section]
\newtheorem{lemma}[theorem]{Lemma}

\newcommand{\qed}{\hfill$\square$\vspace{0.3cm}}

\renewcommand{\theequation}{\arabic{section}.\arabic{equation}}

\title{\Large \bf The Hele-Shaw asymptotics for mechanical \\ models of  tumor growth}

\author{
Beno\^ \i t Perthame\thanks{1- UPMC Paris06, CNRS UMR 7598,
Laboratoire Jacques-Louis Lions, F-75005, Paris. 2- INRIA
Paris-Rocquencourt. Email: benoit.perthame@upmc.fr} \and Fernando
Quir\'os\thanks{Departamento de Matem\'aticas, Universidad
Aut\'onoma de Madrid,  28049-Madrid, Spain.  Email:
fernando.quiros@uam.es} \and Juan Luis V\'{a}zquez\thanks{
Departamento de Matem\'aticas, Universidad Aut\'onoma de Madrid,
28049-Madrid, Spain. Email: juanluis.vazquez@uam.es} }

\date{\today}

\begin{document}

\maketitle


\begin{abstract}
Models of tumor growth,  now commonly used, present several levels
of complexity, both in terms of the biomedical ingredients and the
mathematical description. The simplest ones contain competition for
space using  purely fluid mechanical concepts. Another possible
ingredient is the supply of nutrients through vasculature. The
models can describe the tissue either at the level of cell
densities, or at the scale of the solid tumor, in this latter
case by means of a free boundary problem.

Our first goal here is to formulate a free boundary model of
Hele-Shaw type, a variant including growth terms,  starting from the
description at the cell level and passing to a certain limit. A
detailed mathematical analysis of this purely mechanical model is
performed. Indeed, we are able to prove strong convergence in
passing to the limit, with various uniform gradient estimates; we
also prove uniqueness for the asymptotic Hele-Shaw type problem. The
main tools are nonlinear  regularizing effects for certain porous
medium type equations, regularization techniques \`a la Steklov,
and a Hilbert duality method for uniqueness. At variance with the
classical Hele-Shaw problem, here the geometric motion governed by
the pressure is not sufficient to completely describe the dynamics.
A complete description requires the equation on the cell number
density.

Using this theory as a basis, we go on to consider the more complex
model including nutrients.  We obtain the equation for the limit of
the coupled system; the method relies on some
BV bounds and space/time a priori estimates. Here, new technical
difficulties appear, and they reduce the generality of the results
in terms of the initial data. Finally, we prove uniqueness for the system,
a main mathematical difficulty.

\end{abstract}



\section{Motivation and tumor growth models}
\setcounter{equation}{0}

In the understanding of cancer development, mathematical modeling
and numerical simulations have  nowadays complemented experimental
and clinical observations. The field is now mature; books and
surveys are available, as for example
\cite{Bellomo1,Bellomo2,friedman,Lowengrub_survey}. A first class of
models, initiated in the 70's by Greenspan \cite{greenspan}, has considered
that cancerous cells multiplication is limited by nutrients
(glucosis, oxygen) brought by blood vessels. Models of this class rely on two kinds
of descriptions; either they describe the dynamics of cell
population density \cite{byrne-chaplain} or they consider the \lq
geometric' motion of the tumor through a free boundary problem; see
\cite{{cui_escher}, cui_friedman, friedman_hu} and the references
therein. This stage lasts until the tumor reaches the size of
$\approx 1$mm; then lack of nutrients leads to cell necrosis which
triggers neovasculatures development \cite{chaplain} that supply the
tumor with enough nourishment. This has motivated a new generation of
models where growth is limited by the competition for space
\cite{Bru}, turning the modeling effort towards
mechanical concepts, considering tissues as multiphasic fluids (the
phases could be intersticial water, healthy and tumor cells,
extra-cellular matrix \dots)  \cite{byrne-drasdo, byrne-preziosi,
preziosi_tosin}. This point of view is now sustained by experimental
evidence \cite{JJP}. The term \lq homeostatic pressure', coined
recently, denotes the lower pressure that prevents cell
multiplication by contact inhibition.

The aim of this paper is to explain how asymptotic analysis can link
the two main  approaches, cell density models and free boundary
models, in the context of fluid mechanics. We depart from the
simplest cell population density model, proposed
in~\cite{byrne-drasdo}, in which the cell population density
$\varrho(x,t)$  evolves under pressure forces and cell
multiplication according to the equation
\begin{equation}
\label{model1}
\partial_t\varrho - {\rm div}(\varrho \nabla p) = \varrho \;  \Phi(p),
\end{equation}
where $p$ is the pressure field.
Pressure-limited growth is described by the term $\Phi(p)$,
which  typically satisfies
\begin{equation}
\label{eq:growth}
\Phi'(p) <0 \qquad  \text{and} \qquad \Phi(p_M) =0
\end{equation}
for some $p_M>0$ (the homeostatic pressure). The pressure is assumed to be
a given increasing function of the density.   A representative example is
\begin{equation}
\label{eq:def.pressure}
p=P_m(\varrho):= \frac{m}{m-1}  \left( \frac{\varrho}{\varrho_c}\right)^{m-1}
\end{equation}
with parameter $m >1$. In the free boundary problem to be discussed later, the  value
$\varrho_c$ represents the maximum packing density of cells, as
discussed in \cite{TVCVDP}.

Nutrients consumed by the tumor cells and brought by the  capillary
blood network are a  usual additional ingredient to the modeling.
The situation is described in that case by the system
\begin{equation} \label{model2}
\begin{cases}
\displaystyle\partial_t \varrho -  {\rm div}(\varrho \nabla p) =
\varrho \; \Phi(p,c) ,
\\[2mm]
\displaystyle \partial_t c - \Delta c  = - \varrho \; \Psi(p,c),
\\[2mm]
c(x,t) = c_B >0 \qquad \text{as }\; |x| \to \infty,
\end{cases}
\end{equation}
where $c$ denotes the density of nutrients, and $c_B$ the far field
supply of nutrients (from blood vessels). The coupling functions
$\Phi$,  $\Psi$ are assumed to be smooth and to satisfy the natural
hypotheses
\begin{equation}\label{eq:assumptions.growth.nutrients}
\begin{array}{lll}
\displaystyle\partial_p \Phi< 0,\qquad
&\displaystyle\partial_c\Phi \geq 0, \qquad& \displaystyle \Phi
(p_M, c_B) =0,
\\[10pt]
\displaystyle\partial_p\Psi\leq 0,\qquad
&\displaystyle\partial_c\Psi \geq 0, \qquad&\displaystyle \Psi(p,0)
= 0.
\end{array}
\end{equation}

Variants are possible; for instance, we could assume that nutrients
are  released continuously from a vasculature, then leading to an
equation as
$$
\partial_t c - \Delta c = - \varrho  \; \Psi(p,c) + r
(c_B-c).
$$

We will consider below the purely mechanical model \eqref{model1}
under assumptions \eqref{eq:growth}--\eqref{eq:def.pressure}, and
also the system \eqref{model2}, where nutrients are taken into
account, under assumptions \eqref{eq:def.pressure} and
\eqref{eq:assumptions.growth.nutrients}. In both cases we will show
that the asymptotic limit $m\to\infty$ yields a free boundary model
of Hele-Shaw type, as we were looking for.

Let us recall that  the mathematical theory of the limit for an
equation  like \eqref{model1} in the absence of a growth term is
well developed, the asymptotic limit being in this case the standard
Hele-Shaw model for incompressible fluids with free boundaries.
Early papers on the subject \cite{BBH, Benilan-Crandall, BI1, BI2,
EHKO, FH, FHu, Sa} consider situations in which mass is conserved
and the limit is stationary. In order to have a non-trivial limit
evolution one needs some source, either in the equation or at the
boundary of the domain. A first example, in which there is an
inwards flux at infinity,  is given in~\cite{AGV}, where the authors
study the limit $m\to\infty$ of self-similar focusing (hole-filling)
solutions to the porous medium equation $u_t=\Delta u^m$. A second
example, in which there is a bounded boundary with a nontrivial
boundary data, was first considered in~\cite{GQ1}, and later on
in~\cite{GQ2, I, kim, km}. To our knowledge, the present paper is
the first one in which the evolution in the limit is produced by a
nonlinear source term in the equation. In order to pass to the
limit, three approaches have been used: weak solutions, variational
formulations (using the so-called Baiocchi variable), and viscosity
solutions; see \cite{JK, kim} for this last case. The weak
formulation of Hele-Shaw was first introduced in \cite{DF}, and the
variational formulation in~\cite{EJ}.

Let us also mention that the Hele-Shaw graph, and hence the
Hele-Shaw  equation, can be approximated in other ways, for example
by the Stefan problem. This situation has also been considered in
the literature, even in situations where there is an evolution in
the limit; see for instance~\cite{BKM, GilQuirosVazquez2010, kim}.

We would like to stress that including the growth term is not a
simple change:  several powerful but specific tools do not apply any
longer. More deeply, as we explain later, several approaches have no
chance to work as they are.  This is the reason why we will work on
an equation for the cell density itself rather than the pressure. We
also would like to point out that models of Hele-Shaw type are still
an active field arising in several unexpected applications, see for
instance~\cite{EDS}, and that surface tension is not covered by the
present work; see \cite{CCS,cui_escher, cui_friedman,friedman_hu}.

\medskip

\noindent \textit{Organization of the paper. } We first study the
simplest mechanical model in Section~\ref{sec:mechanical}. The more
complex model with nutrients is studied  in
Section~\ref{sec:nutrients}. The uniqueness proofs for both cases
are highly technical, and are performed separately, in
sections~\ref{sec:uniqueness} and~\ref{sec:nutrient_unique}. These
two sections can safely be skipped by the readers who are mainly
interested in the applications. We finally include an appendix
devoted to some interesting examples for the purely mechanical
model, that illustrate several phenomena.

\medskip

\noindent\textit{Some notations. } We will use several times  the
abridged notations $\varrho(t)$, $p(t)$, meaning
$\varrho(t)(x)=\varrho(x,t)$, $p(t)(x)=p(x,t)$. Given any $T>0$, we
denote $Q_T=\mathbb{R}^N\times(0,T)$, while
$Q=\mathbb{R}^N\times(0,\infty)$.

\section{Purely fluid mechanical model}
\label{sec:mechanical}
\setcounter{equation}{0}

We start our analysis by a detailed investigation of the purely
mechanical model \eqref{model1}, which does not  take into account the consumption of
nutrients, with a pressure field given by \eqref{eq:def.pressure}. A
simple change of scale allows us to assume without loss of
generality  that $\varrho_c=1$. We arrive to the porous medium type
equation, set on $Q$,
\begin{equation}
\label{eq:main} \partial_t\varrho=\Delta
\varrho^m+\varrho\Phi(p),\qquad p=P_m(\varrho):=
\frac{m}{m-1}\varrho^{m-1},\qquad \varrho(0)=\varrho_m^0.
\end{equation}
We assume that  the initial data $\varrho_m^0$ are such that, for
some  $\varrho^0\in L^1_+(\mathbb{R}^N)$,
\begin{equation}
\label{eq:assumptions.initial.data}
\left\{\begin{array}{ll}
\varrho_m^0\ge0,\quad& P_m(\varrho_m^0)\le p_M,\\[10pt]
\|\varrho_m^0- \varrho^0 \|_{L^1(\mathbb{R}^N)}\underset{m \to \infty }{\longrightarrow}  0,\quad& \|\partial_{x_i}
\varrho^0_m\|_{L^1(\mathbb{R}^N)}\le C, \quad i=1,\dots,N.
\end{array}\right.
\end{equation}

\subsection{Main results}

\noindent\textsc{Free boundary limit. } Let $(\varrho_m,p_m)$ be the
unique bounded weak solution to \eqref{eq:main}. We will prove that, along some subsequence, there is a
limit as $m\to\infty$  which turns out to be a solution to a
free boundary problem of Hele-Shaw type.
\begin{theorem}
\label{th:fb} Let $\Phi$ and  $\{\varrho_m^0\}$
satisfy~\eqref{eq:growth} and~\eqref{eq:assumptions.initial.data}
respectively.  Then, after extraction of subsequences, both the
density $\varrho_m$ and  the pressure $p_m$ converge for all $T>0$  strongly in
$L^1(Q_T)$ as $m\to\infty$ to limits $\varrho_\infty \in
C\big([0,\infty); L^1(\mathbb{R}^N)\big)\cap BV(Q_T)$, $p_\infty\in BV(Q_T)$,
that satisfy  $0\leq
\varrho_\infty\leq 1$, $0 \leq p_\infty \leq p_M$, and
\begin{equation}
\label{eq:HS}
\partial_t \varrho_\infty =\Delta
p_\infty +\varrho_\infty  \Phi(p_\infty) \quad\text{in }\mathcal{D}'(Q), \qquad \varrho_\infty(0)
=\varrho^0 \quad\text{in }L^1(\mathbb{R}^N),
\end{equation}
plus the relation $p_\infty\in P_\infty(\varrho_\infty)$, where
$P_\infty$ is the Hele-Shaw monotone graph
\begin{equation}
\label{eq:HS.graph}
P_\infty(\varrho)=\left\{
\begin{array}{ll}
0,\qquad&0\le \varrho<1,\\[3pt]
[0,\infty),\qquad & \varrho=1.
\end{array}
\right.
\end{equation}
\end{theorem}
%
Note that  \eqref{eq:HS.graph} means that a.e.~$P_\infty\ge 0$ and  $P_\infty=0$ a.e.~in  $\{ 0\le
\varrho < 1\}$.

We will obtain also  several important  qualitative properties for
the limit. On the one hand, the \lq tumor is growing' and the pressure
increases, that is,
\begin{equation} \label{eq:rho_tpositive}
\partial_t \varrho_\infty  \geq 0, \qquad
\partial_t p_\infty \geq 0  \qquad \text{in
}\mathcal{D}'(Q).
\end{equation}
On the other hand, if the initial data $\{\varrho_m^0\}$ have a
common compact support, then the limit solution propagates with a
finite speed: $\varrho_\infty(t)$ and $p_\infty(t)$ are compactly
supported for all $t>0$. Finally,  the constructed limit solution enjoys
the monotonicity property, inherited from the case where $m$ is finite,
\begin{equation*} \label{eq:rho_monotone}
\varrho^0 \geq \bar \varrho^0 \; \Longrightarrow \; \left\{
\begin{array}{l}\varrho_\infty(t) \geq \bar \varrho_\infty(t),\\[8pt]
p_\infty(t) \geq \bar p_\infty(t).
\end{array}
\right.
\end{equation*}

\medskip

\noindent\textsc{Transport equation and complementarity formula. }
We will also obtain $L^2$~estimates on the gradients $\nabla p_m$
that will show on the one hand that
$$
\nabla p_\infty\in L^2(Q_T)\quad\text{for all }T>0,
$$
and on the other hand that
$(\varrho_\infty,p_\infty)$ solves a transport
equation.
\begin{theorem}
\label{th:transport} The limit solution $(\varrho_\infty,p_\infty)$
obtained  in Theorem~\ref{th:fb} satisfies
\begin{equation*}
\label{eq:HStransport}
\partial_t \varrho_\infty- {\rm div} \left(\varrho_\infty
\nabla p_\infty \right) = \varrho_\infty \Phi(p_\infty)
\end{equation*}
in a weak sense.
\end{theorem}

A direct calculation shows that the pressure $p_m$ satisfies
\begin{equation}
\label{eq:pressure}
\partial_t p_m=(m-1) p_m \Delta p_m+|\nabla p_m|^2
+(m-1)p_m\Phi(p_m).
\end{equation}
Hence, if we let $m\to\infty$, we formally obtain the \emph{complementarity formula}
\begin{equation}
\label{eq:complementarity.formula}
p_\infty\big(\Delta  p_\infty+\Phi(p_\infty)\big)=0.
\end{equation}
This formal computation can be made rigorous.

\begin{theorem}
\label{thm:complementarity} Under the assumptions of Theorem
\ref{th:fb},   the limit pressure $p_\infty$ satisfies the
complementarity formula
\begin{equation}
\label{eq:complementarity.H1.sense}
\int_{\mathbb{R}^N}\left(-|\nabla
p_\infty|^2+p_\infty\Phi(p_\infty)\right)=0\quad \text{for almost
every } t>0.
\end{equation}
\end{theorem}
It is worth noticing that  \eqref{eq:complementarity.H1.sense} is
equivalent to the strong convergence  of $\nabla p_m$ in $L^2(Q_T)$ for all $T>0$;  see Lemma~\ref{lemma:equivalence.complementarity.strong.convergence}.

\begin{figure}[h!]
 \includegraphics[width=7cm]{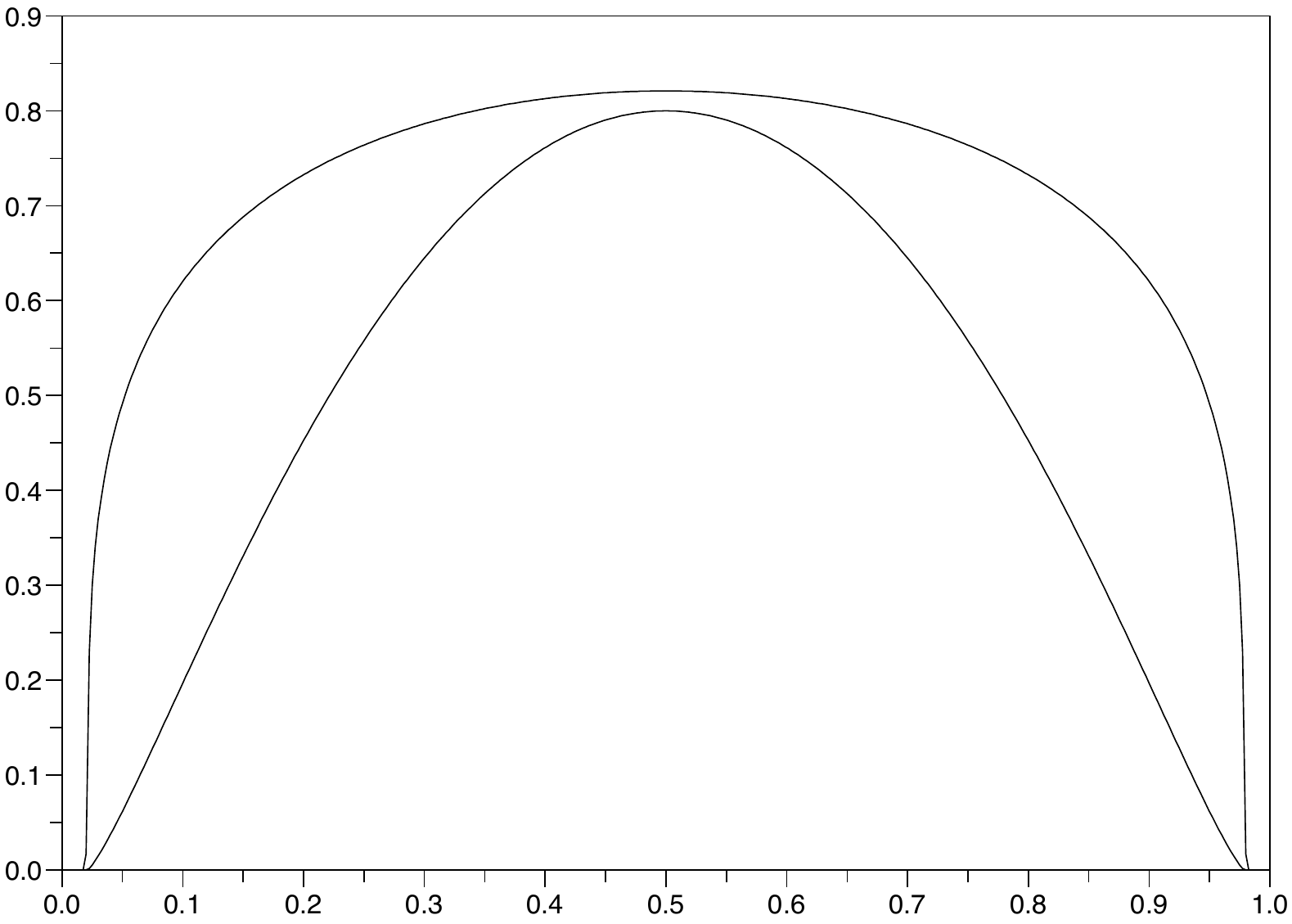} \qquad  \includegraphics[width=7cm]{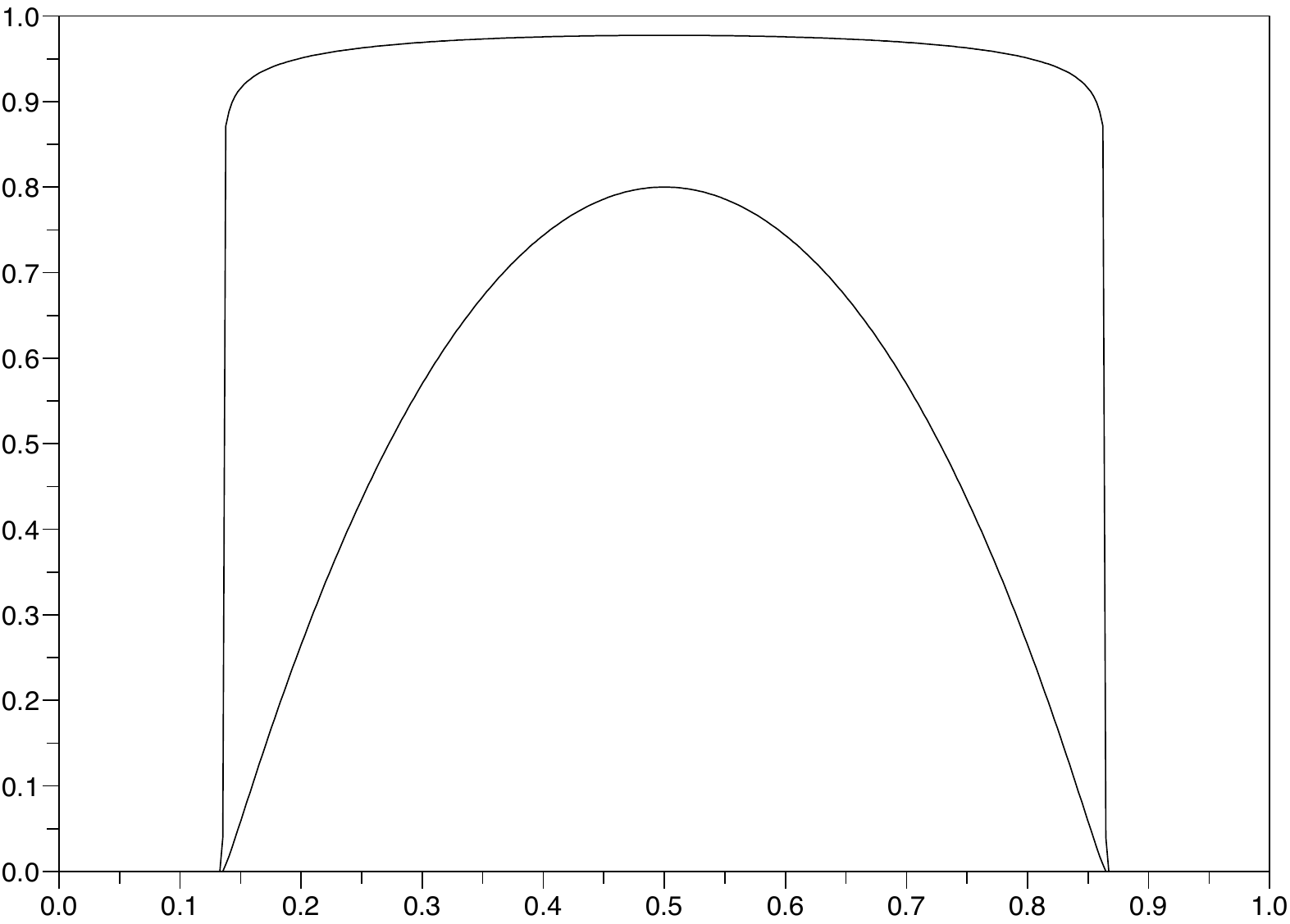}
\\[-25pt]
\caption{\emph{Effect of $m$ large. } A solution to the mechanical
model \eqref{model1},  \eqref{eq:def.pressure} in one dimension with
$\Phi(p)=5(1-p)$. Left: $m=5$. Right: $m=40$. The upper line is
$\varrho$; the bottom line is $p$ (scale enlarged for visibility).
Notice that the density scales are not the same in the two figures.
The initial data is taken with compact support and the solution is
displayed for a time large enough (see Figure~\ref{fig:strict+}
below for and intermediate regime). } \label{fig:strict}
\end{figure}

Since the limit pressure $p_\infty$ is expected to be continuous in
space for all positive times,  the positivity set
$$
\Omega(t):= \{ x; \; p_\infty (x,t) >0 \},
$$
should be well defined for all $t>0$. Notice that it coincides
almost everywhere with the set where $\varrho_\infty =1$; see Figure
\ref{fig:strict}. Indeed, on the one hand, by the definition of the
graph $P_\infty$ we have $\Omega(t)\subset\{x; \; \varrho_\infty
(x,t)=1 \}$; on the other hand, if we had $p_\infty=0$ and $\varrho_\infty =1$ in some set
with positive measure,  then
$\varrho_\infty$ would continue to grow (exponentially) there, which is a
contradiction. Therefore, $\Omega(t)$ may be regarded as the
\emph{tumor}, while the regions where  $0<\varrho_\infty<1$ (mushy
regions, in the literature of phase-changes) correspond to
\emph{precancer} cells.

The complementarity formula~\eqref{eq:complementarity.formula}
indicates that the limit pressure at time $t$, $p_\infty(t)$, should
solve the elliptic equation
\begin{equation}
\label{eq:elliptic.equation} - \Delta p_\infty(t)=
\Phi(p_\infty(t))\ \text{ in }\Omega(t), \qquad p_\infty(t) \in
H^1_0\big(\Omega(t)\big),
\end{equation}
a problem which is wellposed if $\Omega(t)$ is smooth enough, since
$p\mapsto \Phi(p)$  is decreasing. This implies in particular that
in general  regularity in time for the pressure is missing, since time discontinuities may
show up  when two tumors meet; see Subsection~\ref{sec:timedisc} in the Appendix.

\medskip

\noindent\textsc{Geometric motion vs.~equation on the cell number density. }
To complete the description of the limit problem we should be able
to trace $\Omega(t)$ starting from its initial position. The
pressure equation~\eqref{eq:pressure} suggests that we should have
$$
\partial_tp_\infty=|\nabla p_\infty|^2 \quad\text{at }
\partial\Omega(t),
$$
which leads to a geometric motion with normal velocity $V$ at the
boundary of $\Omega(t)$ given by
\begin{equation}
\label{eq:geometric.motion.law}
 V=|\nabla p_\infty|.
\end{equation}
Thus we have arrived to a geometric Hele-Shaw type problem, which
is the classical one when $\Phi=0$.

The above formal computation is expected to be true if we prescribe
a fixed initial pressure $p_m(0)=p^0$ (which implies that the
initial densities converge to the indicator function of the
positivity set of $p^0$). When $\Phi=0$, it was proved to be true in
a viscosity sense in \cite{kim}; see also~\cite{GQ1} for an earlier
result in this direction using a variational formulation of Baiocchi type.

However, if the initial densities $\varrho_m^0$ are such that
$\varrho^0$ is below 1 in a set with positive measure, the result is
no longer true. The main point is that the tumor may meet precancer
zones. At a meeting point $(\overline x, \overline t)$, the example
in Subsection~\ref{subsect:effect.of.eq.density} in the Appendix suggests that the tumor grows faster (also for the case $\Phi=0$),
with a normal velocity given by a rule of the form
\begin{equation}
\label{eq:modified.geometric.motion.law}
 V=\frac{|\nabla p_\infty|}{1-\overline\varrho},
\end{equation}
where $\overline\varrho$ is some limit of $\varrho$ as
$t\to\overline t$ and $x\to\overline x$ from the outside of
$\Omega(t)$. We leave open  this problem, which seems to be a
challenging extension of the viscosity method in \cite{kim}.

Even if we consider the modified geometric motion law
\eqref{eq:modified.geometric.motion.law} instead of
\eqref{eq:geometric.motion.law},  the geometric formulation does not
carry all the information of the limit solution. Indeed, the density
in precancer zones evolves, with an exponential growth, until it
reaches the level $\varrho_\infty=1$, a fact that is not captured
neither by~\eqref{eq:elliptic.equation}, nor
by~\eqref{eq:modified.geometric.motion.law}; see
Figure~\ref{fig:strict+} and the examples in the Appendix.
When $\Phi=0$ we also need~\eqref{eq:HS} to give a full description
of the limit. However, in this latter case the evolution in the mushy
regions is simpler, since the density does not change until it
interacts with the positivity set for the pressure.

\medskip

\begin{figure}[h!]
\includegraphics[width=36mm]{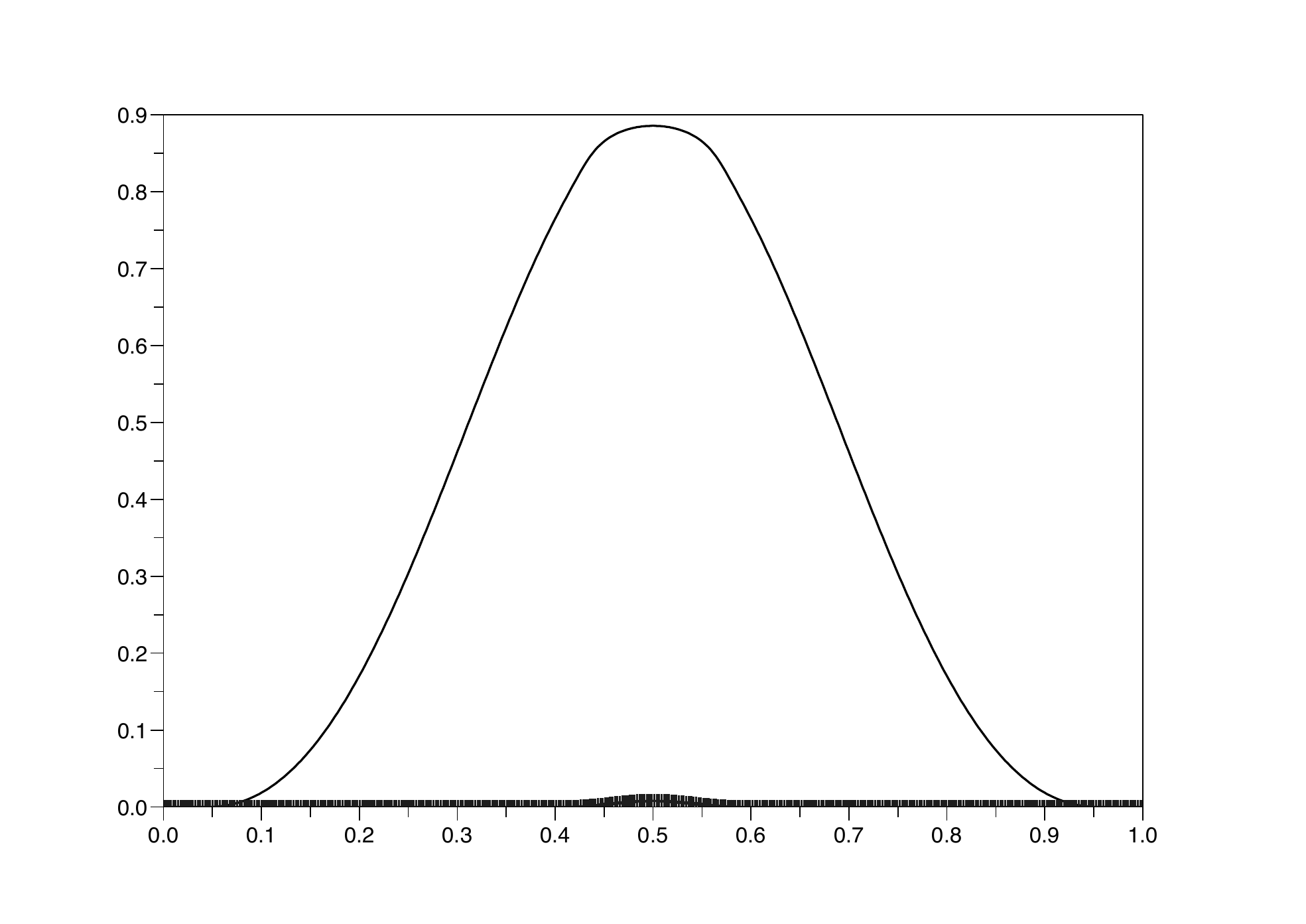}
\includegraphics[width=36mm]{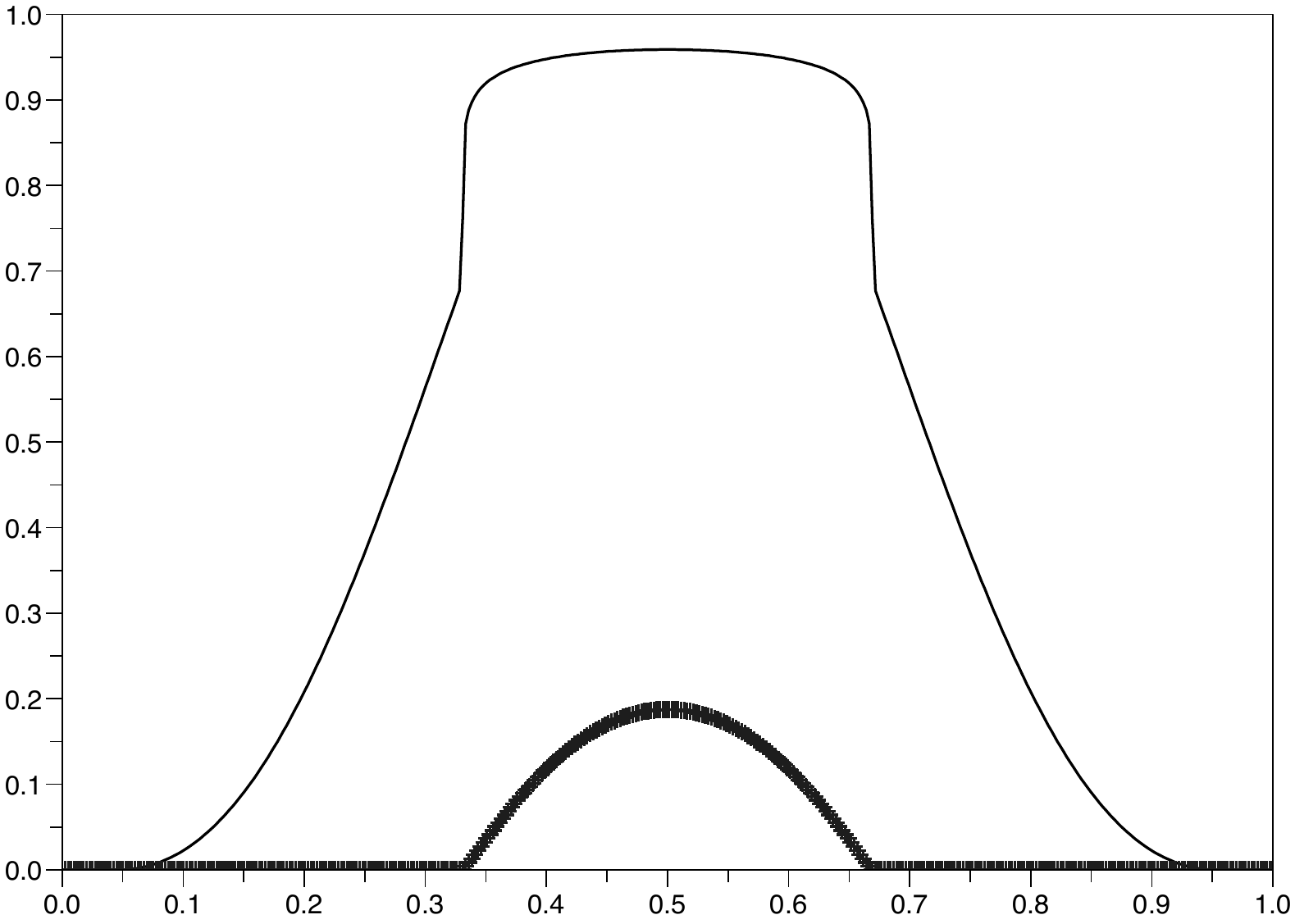}
\includegraphics[width=36mm]{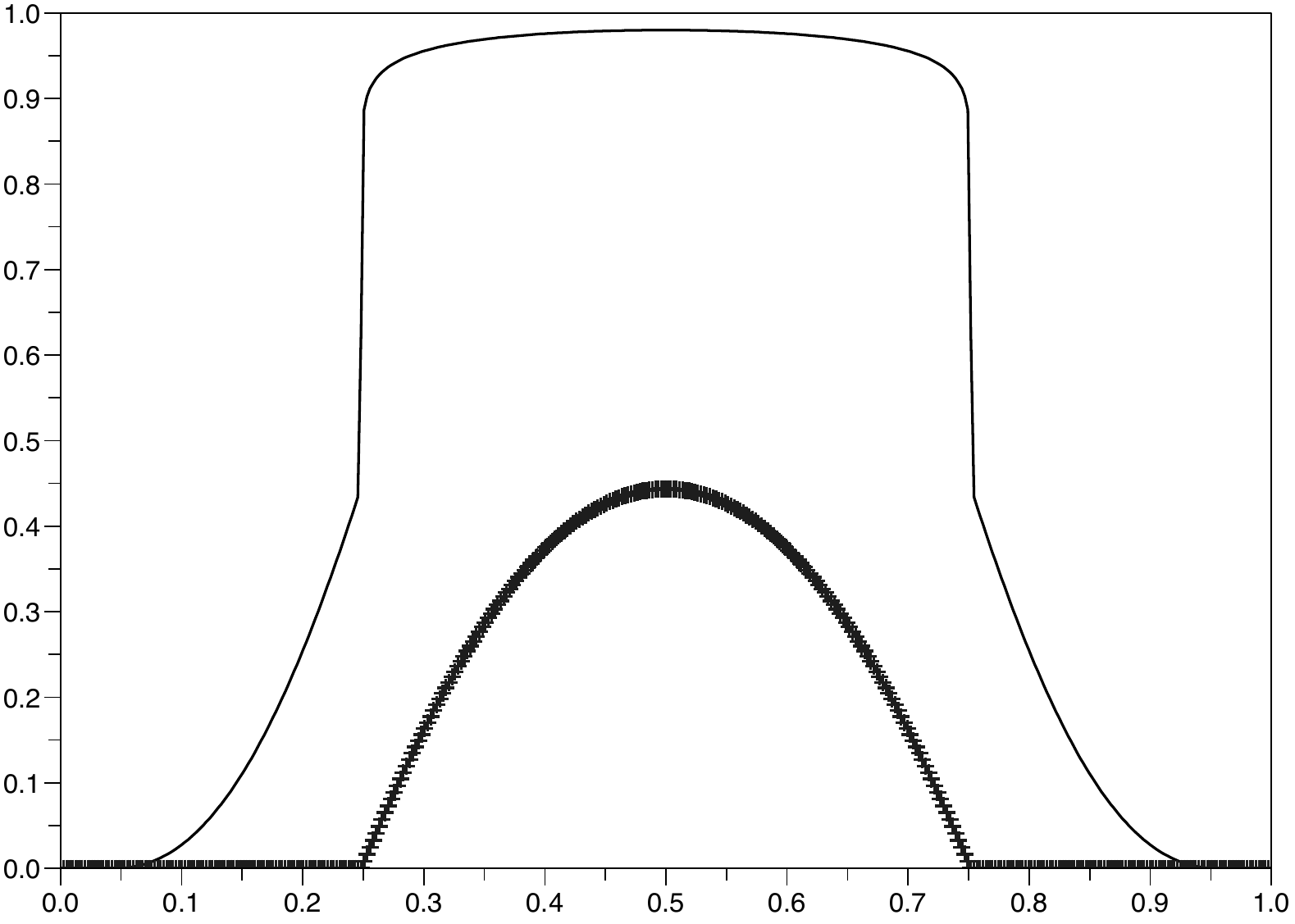}
\includegraphics[width=36mm]{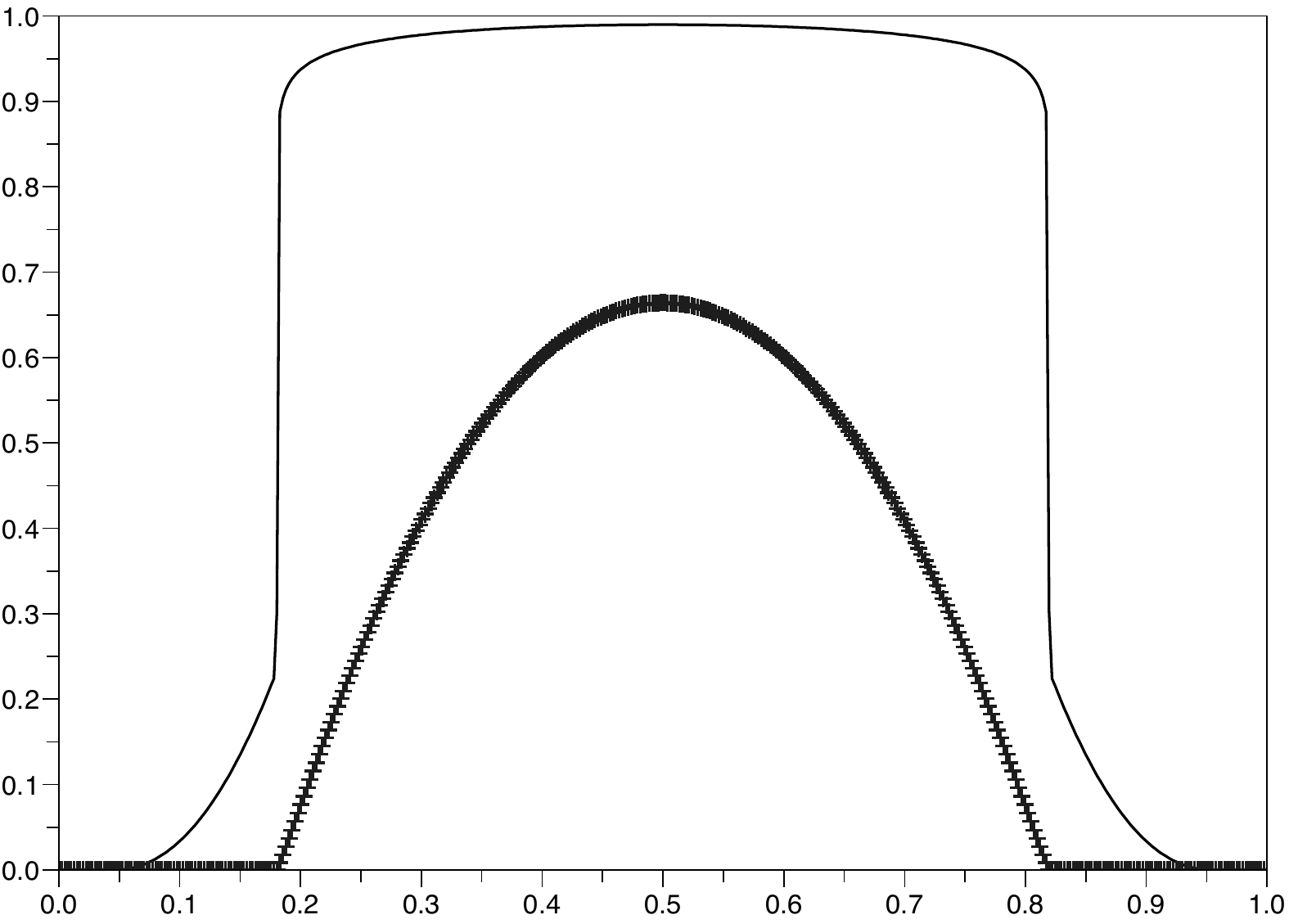}
\\[-25pt]
\caption{\emph{Cell density and pressure carry different
informations. }   Here $m=40$  and the initial data $\varrho$ is
less than $1$. The solution is displayed at four different times. It
shows how the smooth part of $\varrho$ strictly less than $1$ is
growing with $p=0$ (figure on the left). When  $\varrho$  reaches
the value $1$, the pressure becomes positive, increases and creates
a moving  front that delimitates the growing domain where $\varrho
\approx 1$. Thin line is $\varrho$ and thick line is $p$ as functions of $x$. See also
Figure~\ref{fig:tw} in the Appendix.} \label{fig:strict+}
\end{figure}

Let us point out that, in the case of a compactly supported initial data,
precancer zones disappear in finite time.

\medskip

\noindent\textsc{Uniqueness. } To complete the description of  the
asymptotic limit, we conclude with a uniqueness result, relying on a
Hilbert  duality method, for the free
boundary problem~\eqref{eq:HS}.
\begin{theorem}
\label{th:unique} There is a unique pair
$(\varrho,p)$ of functions  in $ L^\infty \big(
(0,\infty);L^1(\mathbb{R}^N)\cap L^\infty ( \mathbb{R}^N)\big) $, $\varrho \in
C\big([0,\infty); L^1(\mathbb{R}^N)\big)$,  $p\in
P_\infty(\varrho)$,  satisfying~\eqref{eq:HS} and such that for all $T>0$:
\begin{itemize}
\item[\rm (i)] $\rho(t)$ is uniformly compactly  supported for $t\in[0,T]$;
\item[\rm (ii)] $|\nabla p| \in L^2 \big( Q_T\big)$;
\item[\rm (iii)] $\partial_t p \in M^1 \big(Q_T\big)$, $\partial_t \varrho \in M^1 \big(Q_T\big)$.
\end{itemize}
\end{theorem}

It is important to notice that if, in addition
to~\eqref{eq:assumptions.initial.data}, we assume that the initial
data $\varrho^0_m$  are uniformly compactly supported, then the limit solution
$(\varrho_\infty,p_\infty)$ to~\eqref{eq:HS} given by
Theorem~\ref{th:fb} falls within the uniqueness class. As a consequence,
convergence is not restricted to a subsequence.

\medskip

\noindent\emph{Remark. } The compact support assumption in
Theorem~\ref{th:unique}  can be removed, to the cost of using rather
technical arguments, in the style of the ones employed in the
uniqueness proof for the case with nutrients,
Theorem~\ref{thm:uniqueness.system} below. We keep it for the sake
of clarity.

\medskip

The rest of this section is devoted to the proofs of the above
mentioned  results, except uniqueness, which is postponed to Section~\ref{sec:uniqueness}.

\subsection{Estimates, strong limit and existence}
\label{sec:estimates}

We consider the solutions to \eqref{eq:main} and proceed to obtain
the estimates that allow us to pass to the limit and prove  Theorem
\ref{th:fb}. This is a standard porous media type equation and all
the manipulations below can be justified, see the monograph \cite{JLV}. The main
observation here is that the lower estimate on $\Delta p_m$ that
does the job when passing to the limit in porous media equations
without a source,  does not work in the present case. However, a
control can be obtained on the quantities  $\Delta p_m +\Phi(p_m)$,
which will be enough for our purposes.

\medskip

\noindent \textsc{$L^\infty$ bounds for $\varrho_m$, $p_m$.  }
Standard comparison arguments yield
\begin{equation*}
\label{eq:L-infty.estimates}
0 \leq \varrho_m\leq \left(\frac{m-1}{m}p_M\right)^{1/(m-1)}
\underset{m \to \infty }{\longrightarrow} 1,\qquad
0\le p_m=P_m(\varrho_m)\le p_M.
\end{equation*}
This  allows in particular to avoid initial layers (which are present when
$\|\varrho^0\|_\infty>1$;  see~\cite{GQ2}).

\medskip

\noindent \textsc{$L^1$ bounds for $\varrho_m$, $p_m$. }
Let $\varrho_m$, $\hat \varrho_m$ be two (non-negative) solutions to
\eqref{eq:main}, and $p_m=P_m(\varrho_m)$, $\hat p_m=P_m(\hat \varrho_m)$ the
corresponding pressures.  We have
$$
\begin{array}{rcl}
\displaystyle\frac{d}{dt}\int_{\mathbb{R}^N}\{
\varrho_m(t)-\hat\varrho_m(t)\}_+
&\le&
\displaystyle \int_{\mathbb{R}^N}\Phi(p_m(t))\{\varrho_m(t)-\hat \varrho_m(t)\}_+\\
&&+\displaystyle
\int_{\mathbb{R}^N}\hat\varrho_m(t)(\Phi(p_m(t))-\Phi(\hat p_m(t)))\;\text{sign}_+(\varrho_m(t)-\hat\varrho_m(t))\\
&\le&\displaystyle\Phi(0)\int_{\mathbb{R}^N}\{\varrho_m(t)-\hat
\varrho_m(t)\}_+,
\end{array}
$$
from where we get
\begin{equation*}
\label{eq:contractivity.scalar.problem}
\int_{\mathbb{R}^N} \{
\varrho_m(t)-\hat\varrho_m(t)\}_+\le
\textrm{e}^{\Phi(0)t}\int_{\mathbb{R}^N} \{
\varrho_m(0)-\hat\varrho_m(0)\}_+.
\end{equation*}
This \lq\lq almost contraction''
property  yields   the
uniform (in $m$) bound
\begin{equation}
\label{eq:mass.estimate} \|\varrho_m(t)\|_{L^1(\mathbb{R}^N)}\le
\textrm{e}^{\Phi(0)t}\|\varrho_m^0\|_{L^1(\mathbb{R}^N)}\le C \textrm{e}^{\Phi(0)t}.
\end{equation}
On the other hand,  $p_m = \frac{m}{m-1} \varrho_m ( \frac{m-1}{m}
p_m)^{\frac{m-2}{m-1}}$. Thus, using \eqref{eq:mass.estimate} we
conclude
\begin{equation*}
\label{eq:pl1bound} \|p_m(t)\|_{L^1(\mathbb{R}^N)} \leq  C
\textrm{e}^{\Phi(0)t} \qquad \text{for } m \ge 2.
\end{equation*}

\medskip

\noindent \textsc{A semiconvexity estimate for $p_m$. } According to the hypotheses~\eqref{eq:growth} on the growth function,
$\displaystyle
r_\Phi=\min_{p\in [0,p_M]}\left(\Phi (p)- p\Phi'(p)\right) >0$.
We will prove that
\begin{equation}
\label{eq:semiconvexity} \Delta p_m(t)
+\Phi(p_m(t))\ge - r_\Phi \textrm{e}^{-(m-1)r_\Phi t}/(1-
\textrm{e}^{-(m-1)r_\Phi t}).
\end{equation}
As a consequence, in the limit we will have
$\Delta p_\infty+\Phi(p_\infty)\ge0$.

We borrow the idea to prove inequality~\eqref{eq:semiconvexity}
from~\cite{AB}, where the case $\Phi=0$ was considered. We write the
equation \eqref{eq:pressure} as
\begin{equation}
\label{eq:p.w}
\partial_t p_m =(m-1)p_m w+|\nabla p_m|^2,\quad\text{with }w=\Delta p_m+\Phi(p_m).
\end{equation}
Let us denote $v= \Delta p_m$. Since $\Phi'<0$, we have
$$
\begin{array}{l}
\partial_t v=(m-1) p_m \Delta w+2(m-1)\nabla p_m
\cdot\nabla w+ (m-1)vw+2\nabla p_m\cdot\nabla
v+2\sum_{i,j}(\partial_{x_ix_j}p_m)^2,
\\[10pt]
\partial_t (\Phi(p_m))=\Phi'(p_m) \partial_t p_m=(m-1)p_m\Phi'(p_m)w+ \Phi'(p_m)|\nabla p_m|^2\\[10pt]
\qquad\qquad\ \, \geq (m-1)p_m\Phi'(p_m)w+ 2 \nabla
(\Phi(p_m)) \cdot \nabla p_m,
\end{array}
$$
which gives
\begin{equation*}
\label{eq:wfinal} \partial_t w\ge(m-1)p_m\Delta w+2m\nabla
p_m\cdot\nabla w+(m-1)w^2-  (m-1) \big( \Phi(p_m)- p_m\Phi'(p_m) \big)w.
\end{equation*}
The function  $W(t)= - r_\Phi \textrm{e}^{-(m-1)r_\Phi t}/(1- \textrm{e}^{-(m-1)r_\Phi
t})$ is a subsolution to   this equation,
and~\eqref{eq:semiconvexity} follows.

\medskip

\noindent\emph{Remarks. } (i)  The right hand side of
\eqref{eq:semiconvexity} behaves as $-1/(t(m-1))$ for $t \approx 0$.

\noindent (ii) For $\Phi\approx 0$ we have $r_\Phi\approx0$, and we
recover the well-known result for the case $\Phi=0$, namely $\Delta
p_m (t)\ge-1/(t(m-1))$.

\medskip

\noindent\textsc{Bounds for $\partial_t p_m$, $\partial_t \varrho_m$. } We
combine~\eqref{eq:semiconvexity} with \eqref{eq:p.w} to obtain an
estimate from below for the time derivative of the pressure,
\begin{equation}
\label{eq:estimate.p_t} \partial_t p_m(t) \ge - (m-1)p_m(t)r_\Phi
\frac{\textrm{e}^{-(m-1)r_\Phi t} }{1- \textrm{e}^{-(m-1)r_\Phi t} }\qquad \text{for }
t>0.
\end{equation}
This in turn gives an estimate from below for the time derivative of
the density,
\begin{equation}
\label{eq:estimate.rho_t}
  \partial_t \varrho_m(t)\ge -   \varrho_m(t)r_\Phi
  \frac{\textrm{e}^{-(m-1)r_\Phi t} }{1- \textrm{e}^{-(m-1)r_\Phi t} }.
\end{equation}
The monotonicity inequalities~\eqref{eq:rho_tpositive} for the limit
problem are then obtained just by letting $m\to\infty$ in
\eqref{eq:estimate.rho_t} and \eqref{eq:estimate.p_t}.

We now use
$|\partial_t\varrho_m|=\partial_t\varrho_m+2\{\partial_t\varrho_m\}_-$
to  obtain
$$
\displaystyle\|\partial_t\varrho_m(t)\|_{L^1(\mathbb{R}^N)}=\displaystyle\frac{d}{dt}\int_{\mathbb{R}^N}
\varrho_m(t)+2\int_{\mathbb{R}^N}\{\partial_t\varrho_m(t)\}_-
\le\displaystyle \left(\Phi(0)+\frac{2
r_\Phi \textrm{e}^{-(m-1)r_\Phi t} }{1- \textrm{e}^{-(m-1)r_\Phi t} }\right)
\|\varrho_m(t)\|_{L^1(\mathbb{R}^N)}.
$$
This, together with~\eqref{eq:mass.estimate}, leads to a uniform
bound in space  and time for time intervals of the form
$t\in\left[\frac1{m-1},T\right]$.

An analogous computation shows that
$$
\begin{array}{rcl}
\displaystyle\int_{\frac1{m-1}}^T\int_{\mathbb{R}^N}|\partial_t
p_m|&\le&\displaystyle\int_{\mathbb{R}^N}p_m(T)-\int_{\mathbb{R}^N}p_m\left(\frac{1}{m-1}\right)\\[4mm]
&&\displaystyle+(m-1)r_\Phi \int_{\frac1{m-1}}^T
\left(\frac{\textrm{e}^{-(m-1)r_\Phi t} }{1- \textrm{e}^{-(m-1)r_\Phi
t}}\int_{\mathbb{R}^N} p_m(t)\right)\,dt \le C(T).
\end{array}
$$

\medskip

\noindent\textsc{$L^1$ bounds for $\nabla \varrho_m$, $\nabla p_m$. }
Let $\alpha=\min_{p\in[0,p_M]}|\Phi'(p)|>0$. We consider the equation
for $\partial_{x_i}\varrho_m$, multiply it by
$\mathop{\rm sign}(\partial_{x_i}\varrho_m)=\mathop{\rm
sign}(\partial_{x_i}p_m)$, and use Kato's inequality;  thanks to the monotonicity of $\Phi$ we obtain
\begin{equation*}
\label{eq:partial.derivative.x}
\begin{array}{rcl}
\partial_t |\partial_{x_i}\varrho_m|-\Delta(m\varrho_m^{m-1}|\partial_{x_i}\varrho_m|)
&\le&\Phi(p_m) |\partial_{x_i}\varrho_m|+\Phi'(p_m)\varrho_m|\partial_{x_i} p_m|\\[4mm]
&\le& \Phi(0)|\partial_{x_i}\varrho_m|-\alpha \varrho_m |\partial_{x_i} p_m|.
\end{array}
\end{equation*}
Integrating in $Q_t$, we get
$$
\|\partial_{x_i}\varrho_m (t)\|_{L^1(\mathbb{R}^N)}+\alpha
\iint_{Q_t}\varrho_m|\partial_{x_i} p_m|\le \|\partial_{x_i}\varrho^0_m \|_{L^1(\mathbb{R}^N)} \textrm{e}^{\Phi(0)t}
\le K\textrm{e}^{\Phi(0)t},
$$
which yields, on the one hand, that
$$
\| \partial_{x_i} \varrho_m(t) \|_{L^1(\mathbb{R}^N)} \leq K
\textrm{e}^{\Phi(0) t},
$$
and on the other hand that
$$
\|\partial_{x_i} p_m\|_{L^1(Q_T)}\le\iint_{Q_T\cap\{\varrho_m\le1/2\}}m\varrho_m^{m-2}|\partial_{x_i}\varrho_m|+
2\iint_{Q_T\cap\{\varrho_m\ge1/2\}}\varrho_m|\partial_{x_i}
p_m|\le C(T).
$$

\medskip

\noindent \textsc{Convergence and identification of the limit. }
Since the families $\varrho_m$ and $p_m$ are bounded in
$W^{1,1}_{\rm loc}(Q)$, we have strong convergence in $L^1_{\rm loc}(Q)$ both for
$\varrho_m$ and $p_m$.

To pass from local convergence to global convergence in $L^1(Q_T)$,
we need  to prove that the mass in an initial strip $t\in[0,1/R]$
and in the tails $|x|>R$ are uniformly (in $m$) small if $R$ is
large enough. The control on the initial strip is immediate using
our uniform, in $m$ and $t$, bounds for
$\|\varrho_m(t)\|_{L^1(\mathbb{R}^N)}$ and
$\|p_m(t)\|_{L^1(\mathbb{R}^N)}$. In order to control the tails, we
consider $\varphi\in C^\infty(\mathbb{R}^N)$ such that
$0\le\varphi\le 1$, $\varphi(x)=0$ for $|x|<R-1$ and $\varphi(x)=1$
for $|x|>R$, and define $\varphi_R(x)=\varphi(x/R)$. Then, for any
$m>2$,
$$
\begin{array}{rcl}
\displaystyle
\frac{d}{dt}\int_{\mathbb{R}^N}\varrho_m(t)\varphi_R&\le&
\displaystyle 2R^{-2}\|\varrho_m(t)\|^{m-1}_{L^\infty(\mathbb{R}^N)}
\|\varrho_m(t)\|_{L^1(\mathbb{R}^N)}\|\Delta\varphi\|_{L^\infty(\mathbb{R}^N)}
+\Phi(0)\int_{\mathbb{R}^N}\varrho_m(t)\varphi_R
\\[4mm]
&\le& \displaystyle CR^{-2}+\Phi(0)\int_{\mathbb{R}^N}\varrho_m(t)\varphi_R.
\end{array}
$$
Therefore, for any $t\in[0,T]$ we have
$$
\begin{array}{rcl}
\displaystyle 0\le
\displaystyle\int_{\mathbb{R}^N}\varrho_m(t)\varphi_R&\le&
\displaystyle
\textrm{e}^{\Phi(0)t}\left(\int_{\mathbb{R}^N}\varrho_m^0\varphi_R+CR^{-2}t
\right)\\[4mm]
\displaystyle &\le&\displaystyle\textrm{e}^{\Phi(0)T}
\left(\|\varrho^0_m-\varrho^0\|_{L^1(\mathbb{R}^N)}+\int_{\mathbb{R}^N}\varrho^0\varphi_R+CR^{-2}T
\right) \le \displaystyle\e
\end{array}
$$
for $R$ and $m$ are large enough. Since
$$
\int_{\mathbb{R}^N}p_m(t)\varphi_R\le  2
\|\varrho_m(t)\|^{m-1}_{L^\infty(\mathbb{R}^N)}
 \int_{\mathbb{R}^N}\varrho_m^0\varphi_R
$$
for any $m\ge 2$, the tail control for the pressures follows easily.

After extraction of
subsequences, we can pass to the a.e.-limit in the equation $
\varrho_mp_m=p_m^{(1+m)/m}$ to obtain that
$(1-\varrho_\infty)p_\infty=0$.  Passing to the limit in the estimates, we have also $0 \le
\varrho_\infty\le 1$, $0\le p_\infty\le p_M$  and
$\varrho_\infty, p_\infty\in BV(Q_T)$ for all $T>0$.

All the above is enough to prove that the pair
$(\varrho_\infty,p_\infty)$  solves the Hele-Shaw equation
\eqref{eq:HS}  but for the question
of the initial data.

\medskip

\noindent \textsc{Time continuity and initial trace. }
We need delicate arguments that are fortunately in the folklore
of the topic of nonlinear diffusion equations, and borrow also
the conclusion from \cite{SV}. Because $\varrho_\infty$ is
non-decreasing in time, we can write for a test function
$0<\zeta(x)<1$ and $0< t_1< t_2 \leq T$,
$$\begin{array}{rl}
\displaystyle\int _{\mathbb{R}^N} |\varrho_\infty(t_2)
-\varrho_\infty(t_1)| \zeta &= \displaystyle\int _{\mathbb{R}^N}
\left[\varrho_\infty(t_2)
-\varrho_\infty(t_1)\right] \zeta
= \displaystyle\int_{t_1}^{t_2} \int_{\mathbb{R}^N} \left[
p_\infty \Delta \zeta + \varrho_\infty \Phi(p_\infty) \zeta\right]
\\[20pt]
&\leq C(T) (t_2-t_1)\left( \| \Delta \zeta \|_\infty +1\right).
\end{array}
$$
Taking a sequence of such uniformly smooth functions  that converge
to $1$, we find that $\varrho_\infty \in C\big([0,\infty);
L^1(\mathbb{R}^N)\big)$ (in fact with a locally uniform Lipschitz
constant).

In order to identify the initial trace, we observe that for any test function $\zeta$ as above,
$$
\int_{\mathbb{R}^N}\varrho_m(t)\zeta-\int_{\mathbb{R}^N}\varrho^0_m\zeta=\int_0^t\int_{\mathbb{R}^N}
\left(p_m\Delta\zeta+\varrho_m\Phi(p_m)\right)\zeta.
$$
Letting $m\to\infty$, we have
$$
\int_{\mathbb{R}^N}\varrho_\infty(t)\zeta-\int_{\mathbb{R}^N}\varrho^0_\infty\zeta=\int_0^t\int_{\mathbb{R}^N}
\left(p_\infty\Delta\zeta+\varrho_\infty\Phi(p_\infty)\right)\zeta.
$$
Letting first $t\to0$ and then $\zeta\to1$, we conclude that $\varrho(0)=\varrho^0$ in $L^1(\mathbb{R}^N)$.

\subsection{$L^2$ estimates for $\nabla p_m$, transport equation and complementarity formula}
\label{sec:obstacle}
We now prove Theorems~\ref{th:transport} and~\ref{thm:complementarity}, both of them connected with the fact that $\nabla p_m\in L^2(Q_T)$ for all $T>0$.

Since we have already proved that both $\varrho_m$ and $p_m$ converge strongly in $L^1(Q_T)$ for all~$T>0$, the first of these theorems just depends on obtaining a uniform bound in $L^2(Q_T)$ for~$\nabla p_m$.

\medskip

\noindent\emph{Proof of Theorem~\ref{th:transport}. }
We rewrite the pressure equation~\eqref{eq:pressure} as
\begin{equation}
\label{eq:pressure.2}
\partial_t p_m =\frac{m-1}2\Delta p_m^2-(m-2)|\nabla p_m|^2+(m-1)p_m\Phi(p_m).
\end{equation}
Integrating in $Q_T$ we obtain the required estimate,
\begin{equation*}
\label{eq:estimate.p_x} \|\nabla p_m\|_{L^2(Q_T)}^2\le\frac{m-1}{m-2}\Phi(0) \|p_m\|_{L^1(Q_T)}
+\frac{1}{m-2}\|p^0_m\|_{L^1(\mathbb{R}^N)}.
\end{equation*}
\qed

\noindent\emph{Remark. } Combining~\eqref{eq:pressure.2}
with~\eqref{eq:estimate.p_t}, we obtain
\begin{equation*}
\label{eq:estimate.p_x1} \|\nabla p_m(t)\|_{L^2(\mathbb{R}^N)}^2\le
\int_{\mathbb{R}^N}p_m(t) \left(\frac{m-1}{m-2}\Phi(0) +(m-1) r_\Phi
  \frac{\textrm{e}^{-(m-1)r_\Phi t} }{1- \textrm{e}^{-(m-1)r_\Phi t} }\right)
\end{equation*}
for all $m>2$. This gives a uniform estimate for $\|\nabla p_m(t)\|_{L^2(\mathbb{R}^N)}$, $t\in\left[\tau,T\right]$,
for any fixed values $0<\tau<T$. However, it does not allow to go down to $t=0$.

\medskip

From the  equation in~\eqref{eq:main} and the $L^1$ bound on
$\partial_t \varrho_m$ we conclude that
\begin{equation*}
\label{eq:estimate.p_xx} \left\| \Delta \varrho_m^m(t) \right\|_{L^1(\mathbb{R}^N)} \leq C
\textrm{e}^{\Phi(0) t}.
\end{equation*}
This is an optimal bound since $p_\infty$ has \lq corners' on
$\partial \Omega(t)$.  Because this only gives space compactness, it
is not enough to establish the  complementarity formula
\eqref{eq:complementarity.formula} on the pressure when passing to
the limit in  \eqref{eq:pressure}. We will need to perform a time
regularization argument {\em \`a la Steklov}.

\medskip

\noindent\emph{Proof of Theorem~\ref{thm:complementarity}. } Let
$\omega_\e(t) \geq 0$ be  a regularizing kernel in time with support $(-\e,0)$.
We use the notation $\varrho_{m,\e}(t) =
\omega_\e\star\varrho_m
(t)=\int_{\mathbb{R}}\omega_\e(t-s)\varrho_m(s)\,ds$. The equation
\eqref{eq:main} gives
\begin{equation}
\label{eq:regularized} \Delta (\varrho_m^m\star \omega_\e) =\partial_t
\varrho_{m,\e}  -( \varrho_m\Phi(p_m) ) \star \omega_\e.
\end{equation}
From here it follows that  $U_m=\Delta( \varrho_m^m\star\omega_\e )$
is uniformly  (in $m$) smooth in time and $H^1$ in space for $\e$
fixed. Therefore $U_m$ converges strongly in $L^1_{\rm loc}(Q)$. Hence, after multiplying \eqref{eq:regularized} by
$p_m$,  we may pass to the limit to obtain
$$
p_\infty \big(\Delta (p_\infty \star \omega_\e) +( \varrho_\infty
\Phi(p_\infty) ) \star \omega_\e\big) = \lim_{m \to \infty } p_m
\partial_t \varrho_{m,\e}.
$$
In order to estimate the right hand side we  make the following
decomposition,
$$
\begin{array}{l}
 (p_m \partial_t \varrho_{m,\e})(t)=
\underbrace{\frac{m}{m-1}\int_{\mathbb{R}}\varrho_m^{m-1}(s)  \partial_t
\varrho_{m}(s) \omega_\e(t-s)ds}_{{\text{I}}_m(t)}\\[5mm]
\qquad\qquad\qquad
 +\underbrace{\frac{m}{m-1}\int_{\mathbb{R}}\left(\varrho_m^{m-1}(t)
-\varrho_m^{m-1}(s)\right)  \left(\partial_t \varrho_{m}(s) + \frac Cs\right)
\omega_\e(t-s)ds}_{{\text{II} }_m(t)}\\[5mm]
\qquad\qquad\qquad\underbrace{-\frac{Cm}{m-1}\int_{\mathbb{R}}
\left(\varrho_m^{m-1}(t) -\varrho_m^{m-1}(s)\right)
\frac{\omega_\e(t-s)}s\,ds}_{{\text{III} }_m(t)},
\end{array}
$$
where $C$ is a constant such that $\partial_t \varrho_{m} + \frac Ct
\geq 0$; see estimate~\eqref{eq:estimate.rho_t}.

For the first term we have
$$
\int_{\mathbb{R}^N}|{\text{I}}_m(t)|\le\frac 1{m-1}\int_{t}^{t+\e}
\omega_\e(t-s)\left(\int_{\mathbb{R}^N}|\partial_t
\varrho_m^m(s)|\right)\, ds\le \frac{C_\e}{m-1}.
$$
Regarding the second term, the estimate~\eqref{eq:estimate.p_t} implies
that $\partial_t p_m\ge -C$. Therefore, $\varrho_m^{m-1}(t)
-\varrho_m^{m-1}(s)\le C\e$.  Let $\zeta\in
\mathcal{D}'(Q)$, and $\tau$  the
smallest time in its support. Then,
$$
\displaystyle\iint_{Q}\zeta \;  {\text{II}}_m   \leq   C\e
\iint_Q\int_{\mathbb{R}}   \zeta \; \left(\partial_t \varrho_{m}(s)+
\frac C\tau \right) \omega_\e(t-s)\,dsdx dt\le C\e.
$$

The third term is easy to treat.  Since $s\ge t>0$, for any test
function  $\zeta$ as above,
$$
\displaystyle\iint_Q\zeta \; {\text{III}}_m \underset{m \to
\infty}{\longrightarrow}  -C  \iint_Q\zeta \int_{\mathbb{R}} \left(
p_\infty(t) -p_\infty(s)\right) \frac{\omega_\e(t-s)}s\,dsdxdt=o(1)
$$
as $\e \to 0$.

These three calculations give, in the limit $m\to \infty$,
$$
p_\infty  \big( \Delta (p_\infty\star \omega_\e) + ( \varrho_\infty
\Phi(p_\infty) )  \star \omega_\e\big) \leq o(1) \quad \text{as }\;
\e \to 0 \quad \text{ in } \; \mathcal{
D}'(Q).
$$
To recover the desired information,
$$
p_\infty  \big(\Delta p_\infty + \Phi(p_\infty)\big) \leq 0,
$$
it remains to pass to the limit as $\e \to 0$ after noticing that
$p_\infty\varrho_\infty = p_\infty $ and that the differential
term can be treated through its weak formulation; indeed, after
testing it against a test function as before, it is written as
$$
-\iint_Q  \big[ \zeta \nabla p_\infty  \cdot \nabla (p_\infty\star \omega_\e)  +
 p_\infty   \nabla \zeta \cdot \nabla (p_\infty\star \omega_\e) \big]
$$
which passes to the limit $\e\to0$ because we already know that
$\nabla p_\infty\in L^2(Q_T)$. Therefore,
$$
-\iint_{Q}  \big[ \zeta |\nabla p_\infty|^2 + p_\infty \nabla \zeta
\cdot \nabla p_\infty+\zeta p_\infty\Phi(p_\infty) \big]=0.
$$
\qed

\medskip

As we have already mentioned, the complementarity formula is
equivalent  to the strong convergence of the gradients.
\begin{lemma}
\label{lemma:equivalence.complementarity.strong.convergence} The
complementarity formula \eqref{eq:complementarity.H1.sense} holds if
and only if $\nabla p_m \underset{m \to \infty}{\longrightarrow}
\nabla p_\infty$  strongly in $L^2(Q_T)$ for all $T>0$,
\end{lemma}

\noindent\emph{Proof. } We consider a test function
$\zeta\in\mathcal{D}'(Q)$  and use
\eqref{eq:p.w} and~\eqref{eq:estimate.p_x} to obtain
$$
\iint_Q  \zeta p_m (\Delta p_m+\Phi(p_m)) = \frac{1}{m-1} \iint_Q
\zeta \left(\partial_t p_m - |\nabla p_m|^2 \right) \underset{m \to
\infty}{\longrightarrow} 0.
$$
This implies, after an integration by parts, that
$$
\iint_Q \left[ -\zeta  |\nabla p_m|^2 - p_m \nabla \zeta\cdot\nabla
p_m +p_m  \zeta \Phi(p_m) \right] \underset{m \to
\infty}{\longrightarrow} 0.
$$
Let $r= \liminf_{m\to\infty}  \iint_Q \zeta|\nabla p_m|^2 - \iint_Q
\zeta |\nabla p_\infty|^2$. Thanks to Fatou's Lemma, $r\ge0$. On the
other hand, since $p_m$ converges strongly, after extraction of
subsequences we have
$$
\iint_Q \left[ -\zeta  |\nabla p_\infty|^2 - p_\infty \nabla
\zeta\cdot \nabla p_\infty + \zeta p_\infty \Phi(p_\infty) \right] =
r,
$$
and $r=0$ is equivalent to the strong convergence.
\qed

\subsection{Finite speed of propagation}
\label{sec:fsp}

If the initial data are compactly supported uniformly in $m$, we can
control the supports of $p_m(t)$ and $\varrho_m(t)$ uniformly in $m$
and $t\in [0,T]$ for all $T>0$. This implies that the speed of
propagation is finite for the limit problem.

The control is performed comparing the pressures with functions of the form
\begin{equation}
\label{eq:def.supersolution}
P(x,t)=\left(C-\frac{|x|^2}{4(\tau+t)}\right)_+.
\end{equation}
The key point  is that these functions, which are viscosity
solutions of the Hamilton-Jacobi equation $P_t=|\nabla P|^2$,  are
supersolutions to the equation
\begin{equation}
\label{eq:pressure.simplified}
\partial_t p=(m-1)p\Delta p+|\nabla p|^2+(m-1)p\Phi(0)
\end{equation}
for some time interval which does not depend on $m$. This idea was
already used for the case $\Phi=0$ in \cite{GQ1,GQ2}. However, in
the absence of reaction, $P$ is a supersolution for all times, and
the proof is a bit easier.

\begin{lemma}
Let $\{\varrho_m^0\}$ satisfying~\eqref{eq:assumptions.initial.data}
and such that all their supports are contained in a common ball
$B_R$. Then, for every $T>0$, there is a radius $\mathcal{R}_T$
depending only on $\sup_{m}\|\varrho_m^0\|_{L^\infty(\mathbb{R}^N)}$, $R$ and $T$,
such that the supports of $\{\varrho_m(t) \}$ are contained in
the ball $B_{\mathcal{R}_T}$ for all $m>1$ and $t\in[0,T]$.
\end{lemma}

\noindent{\it Proof. } We consider $P$ as
in~\eqref{eq:def.supersolution}  with $\tau=N/(4\Phi(0))$ and $C$
large enough so that $P(x,0)\ge \varrho_m^0(x)$ for all
$x\in\mathbb{R}^N$ and $m>1$. An easy computation shows that
$$
\partial_t  P - (m-1)P\Delta P- |\nabla P|^2 - (m-1)P\Phi(0)=(m-1)P(x,t)\left(\frac{N}{2(\tau+t)}-\Phi(0)\right)\ge0
$$
for all $t\in[0,\frac{N}{4\Phi(0)}]$. Since the functions $p_m$ are
subsolutions to~\eqref{eq:pressure.simplified}, we conclude that the
supports of the functions $\varrho_m$ are contained in the ball of
radius $\sqrt{\frac{2CN}{\Phi(0)}}$ for all
$t\in[0,\frac{N}{4\Phi(0)}]$. Moreover, $\|p_m(\frac{N}{4\Phi(0)})\|_{L^\infty(\mathbb{R}^N)}\le C$, and the argument can be iterated
to reach the time $T$ in a finite number of steps. \qed

\section{Uniqueness for the purely mechanical limit problem}
\label{sec:uniqueness}
\setcounter{equation}{0}

In this section we prove, under suitable  assumptions,
that the limit problem~\eqref{eq:HS} has a unique solution, Theorem
\ref{th:unique}. The main difficulty comes from the fact that $p$ is
not a Lipschitz, single-valued function of $\varrho$. Hence, we
cannot apply directly the ideas developed in
\cite{Benilan-Crandall-Pierre-1984} to adapt Hilbert's duality method to the porous medium equation. The
technique has to be modified, in the spirit of \cite{Crowley-1979}.

\medskip

\noindent\textsc{Dual problem. } Consider two
solutions, $(\varrho_1,p_1)$, $(\varrho_2,p_2)$. Let $\Omega$ be a
bounded domain containing the supports of both solutions for all
$t\in[0,T]$, and $\Omega_T=\Omega\times(0,T)$. Then we have
\begin{equation}
\label{eq:uniqueness.distributional.solution}
\iint_{\Omega_T}\big[(\varrho_1-\varrho_2)
\partial_t\psi+(p_1-p_2)\Delta
\psi+(\varrho_1\Phi(p_1)-\varrho_2\Phi(p_2))\psi\big]=0
\end{equation}
for all suitable test functions $\psi$. This can be rewritten as
\begin{equation}
\label{eq:weak.formulation.limit.problem}
\iint_{\Omega_T}(\varrho_1-\varrho_2+p_1-p_2)\left[ A\partial_t\psi+{B}\Delta
\psi+A\Phi(p_1)\psi - C {B} \psi \right]=0,
\end{equation}
where, for some fixed $\nu >0$,
$$
\begin{array}{l}
0 \leq \displaystyle A=\frac{\varrho_1-\varrho_2}{(\varrho_1-\varrho_2)+(p_1-p_2)} \leq 1,\\[10pt]
0 \leq \displaystyle B=\frac{p_1-p_2}{(\varrho_1-\varrho_2)+(p_1-p_2)} \leq  1,\\[10pt]
0 \leq \displaystyle C= \; - \varrho_2\,\frac{\Phi(p_1)-\Phi(p_2)}{p_1-p_2}\;  \leq  \nu.
\end{array}
$$
To arrive to these bounds on $A$ and $B$,  we define $A=0$  when
$\varrho_1 =\varrho_2$, even when $p_1=p_2$, and $B=0$ when
$p_1=p_2$, even when $\varrho_1=\varrho_2$.

The idea of Hilbert's duality method is to solve the \emph{dual problem}
$$
\left\{
\begin{array}{l}
 A\partial_t\psi+{B}\Delta
\psi+A\Phi(p_1)\psi - C {B} \psi=AG\quad \text{in }\Omega_T,\\[10pt]
\psi=0\quad\text{in }\partial\Omega\times(0,T),\quad
\psi(\cdot,T)=0\quad\text{in }\Omega,
\end{array}
\right.
$$
for any smooth function $G$, and use $\psi$ as test function. This would yield
\begin{equation*}
\label{eq:uniqueness.density}
\iint_{\Omega_T}(\varrho_1-\varrho_2)G=0,
\end{equation*}
from where uniqueness for the density is immediate. Uniqueness for
the pressure would then  follow
from~\eqref{eq:uniqueness.distributional.solution}.

However, on the one hand the coefficients of the dual problem are
not smooth, and, on the other, since   $A$ and $B$ are not strictly
positive, the dual equation is not uniformly parabolic. A smoothing
argument  is required.

\medskip

\noindent\textsc{Regularized dual problem. } Let $\{A_n\}$,
$\{B_n\}$, $\{C_n\}$, $\{\Phi_{1,n}\}$ be sequences of smooth
bounded functions such that
\begin{equation*}
\label{eq:choice.approximations}
\left\{\begin{array}{lll}
\|A-A_n\|_{L^2(\Omega_T)}<\alpha/n,\qquad&1/n<A_n\le 1,\\[3mm]
\|B-B_n\|_{L^2(\Omega_T)}<\beta/n,\qquad&1/n<B_n\le 1,\\[3mm]
\|C-C_n\|_{L^2(\Omega_T)}<\delta_1/n,\qquad&0 \leq C_n <M_1,\qquad&\| \partial_t C_n \|_{L^1(\Omega_T)} \leq K_1 \\[3mm]
\|\Phi(p_1)-\Phi_{1,n}\|_{L^2(\Omega_T)}<\delta_2/n,\qquad&|\Phi_{1,n}|<M_2,\qquad&
\| \nabla \Phi_{1,n}  \|_{L^2(\Omega_T)} \leq K_2,
\end{array}
\right.
\end{equation*}
for some constants $\alpha,\beta,\delta_1,\delta_2,M_1,M_2,K_1,K_2$.

For any smooth function $G$, we consider  the solution
$\psi_n$ of the backward heat equation
\begin{equation}
\label{eq:auxiliary.problem}
\left\{
\begin{array}{l}
\displaystyle\partial_t \psi_n+\frac{B_n}{A_n}\Delta \psi_n+\Phi_{1,n}\psi_n - C_n\frac{B_n}{A_n}\psi_n=G \quad \text{in }\Omega_T,\\[10pt]
\psi_n=0\quad\text{in }\partial\Omega\times(0,T),\quad
\psi_n(\cdot,T)=0\quad\text{in }\Omega.
\end{array}
\right.
\end{equation}
The coefficient $B_n/A_n$ is continuous, positive and  bounded below
away from zero. Thus, the equation satisfied by $\psi_n$ is
uniformly parabolic in $\Omega_T$. Hence $\psi_n$ is smooth and can be
used as a test function in  \eqref{eq:weak.formulation.limit.problem}.

Combining \eqref{eq:weak.formulation.limit.problem} and
\eqref{eq:auxiliary.problem}, we have
$$
\displaystyle\iint_{\Omega_T} (\varrho_1-\varrho_2)G
=I^1_n-I^2_n-I^3_n+I^4_n,
$$
where
$$
\begin{array}{l}
I^1_{n}=\displaystyle\iint_{\Omega_T} \big((\varrho_1-\varrho_2)+(p_1-p_2)\big) \frac{B_n}{A_n}(A-A_n)\left(\Delta \psi_n -C_n \psi_n\right),\\[10pt]
I^2_{n}=\displaystyle\iint_{\Omega_T} \big((\varrho_1-\varrho_2)+(p_1-p_2)\big) (B-B_n) \left(\Delta\psi_n - C_n \psi_n\right) ,\\[10pt]
I^3_{n}=\displaystyle\iint_{\Omega_T} (\varrho_1-\varrho_2)(\Phi(p_1)-\Phi_{1,n})\psi_n,\\[10pt]
I^4_{n}=\displaystyle\iint_{\Omega_T}
\big((\varrho_1-\varrho_2)+(p_1-p_2)\big) B(C-C_n)\psi_n.
\end{array}
$$

\medskip

\noindent\textsc{Limit $n\to\infty$. Uniform bounds and conclusion. }
Our aim is to prove that $\lim_{n\to\infty}I^i_{n}=0$, $i=1,\dots,4$,
which implies that $\varrho_1=\varrho_2$. For this task we will use
some bounds that we gather  in the following lemma.

\begin{lemma}
\label{lem:bounds.uniqueness.scalar.equation}
There are constants $\kappa_i$, $i=1,2,3$, depending on $T$ and $G$, but not on $n$, such that
$$
\displaystyle\| \psi_n\|_{L^\infty(\Omega_T)}\le \kappa_1,  \quad \sup_{0\leq t \leq
T} \|\nabla \psi_n (t) \|_{L^2(\Omega)}\le\kappa_2,\quad
\|(B_n/A_n)^{1/2} \left( \Delta \psi_n -C_n
\psi_n \right)\|_{L^2(\Omega_T)} \leq \kappa_3.
$$
\end{lemma}

\noindent\emph{Proof. } The first bound is just the maximum
principle because $C_n$ is non-negative  and $ \Phi_{1,n}$ uniformly
bounded.

The second bound is obtained multiplying the equation
in~\eqref{eq:auxiliary.problem} by $\Delta \psi_n - C_n \psi_n$.
After integration in $ \Omega\times(t,T)$, we get
\begin{equation}
\label{eq:estimates.auxiliary.problem}
\begin{array}{l}
\displaystyle \frac 12 \| \nabla \psi_n (t) \|_{L^2(\Omega)}^2+\int_t^T
\int_\Omega  \frac{B_n}{A_n} | \Delta \psi_n -C_n \psi_n |^2=-\int_{\Omega}\left(\frac{C_n\psi_n^2}{2}\right)(t)
\\[4mm]
\displaystyle \quad +\int_t^T \int_\Omega \left[ -\partial_t C_n
\frac{ \psi_n^2}{2} -  \Phi_{1,n} |\nabla \psi_n|^2 - \psi_n \nabla
\Phi_{1,n}  \cdot\nabla \psi_n+ C_n
\Phi_{1,n} \psi_n^2 + \psi_n \Delta G - G C_n \psi_n  \right]\\[4mm]
\displaystyle\   \leq K\left(1-t+  \int_t^T \|\nabla \psi_n (s) \|_{L^2(\Omega)}^2\, ds\right),
\end{array}
\end{equation}
where we denote by $K$ various constants independent of $n$.
We now use the Gr\"{o}nwall lemma to conclude our second uniform bound.
Then we can use again the inequality \eqref{eq:estimates.auxiliary.problem} to conclude
the third one. \ $\hfil$\qed

We now get
$$
\begin{array}{rcl}
\displaystyle |I^1_{n}|&\le& \displaystyle K\iint_{\Omega_T}
\frac{B_n} {A_n}|A-A_n|\,|\Delta\psi_n - C_n \psi_n |\le
K\|(B_n/A_n)^{1/2}(A-A_n)\|_{L^2(\Omega_T)}\\[10pt]
&\le&\displaystyle
Kn^{1/2}\|A-A_n\|_{L^2(\Omega_T)}\le K\alpha/n^{1/2};\\[10pt]
\displaystyle |I^2_{n}|&\le& \displaystyle K\iint_{\Omega_T}|B-B_n|\,
|\Delta\psi_n - C_n \psi_n|\le
K\|(A_nn^{1/2}/B_n)^{1/2}(B-B_n)\|_{L^2(\Omega_T)}\\[10pt]
&\le&\displaystyle
Kn^{1/2}\|B-B_n\|_2\le K\beta/n^{1/2};\\[10pt]
\displaystyle |I^3_{n}|&\le& \displaystyle\iint_{\Omega_T}
|\varrho_1-\varrho_2|\, |\Phi(p_1)-\Phi_{1,n}|\,|\psi_n| \le
K\|\Phi(p_1)-\Phi_{1,n}\|_{L^2(\Omega_T)}
\le K\delta_2/n;
\\[10pt]
\displaystyle |I^4_{n}|&\le& \displaystyle K\iint_{\Omega_T} B
|C-C_n|\,|\psi_n|  \le \displaystyle K \|C-C_n\|_{L^2(\Omega_T)} \le K/n.
\end{array}
$$

Once we have proved that $\varrho_1=\varrho_2$,
equation~\eqref{eq:uniqueness.distributional.solution} says that
$$
\iint_{\Omega_T} \ \big[ (p_1-p_2)\Delta \psi
+\varrho_1(\Phi(p_1)-\Phi(p_2))\psi \big] = 0,
$$
from where the result will follow by taking as test function
$\psi=p_1-p_2$,  using the monotonicity of $\Phi$.

\section{Model with nutrients}
\label{sec:nutrients}
\setcounter{equation}{0}

When nutrients are taken into account, we are led to problem
\eqref{model2}. We  again assume that the pressure field is given
by~\eqref{eq:def.pressure},  with $\varrho_c$ set equal to 1.

Our assumptions on the initial data  are  stronger than for the
purely mechanical problem. Namely, in addition to
\eqref{eq:assumptions.initial.data}, we assume that for some $c^0$
such that $c_B-c^0\in L^1_+(\mathbb{R}^N)$,
\begin{equation}\label{eq:assumptions.initial.data.nutrients}
\left\{
\begin{array}{l}
0\le c_m^0 < c_B, \quad
\|c_m^0-c^0\|_{L^1(\mathbb{R}^N)}\underset{m \to \infty
}{\longrightarrow}  0,\quad \|(c_m^0)_{x_i}\|_{L^1(\mathbb{R}^N)}\le
C,\quad
i=1,\dots,N,\\[8pt]
\|{\rm div}(\varrho_m^0 \nabla p_m^0) + \varrho_m^0 \;
\Phi(p_m^0,c_m^0)\|_{L^1(\mathbb{R}^N)}\le C, \quad  \|\Delta c_m^0 -
\varrho_m^0  \; \Psi(p_m^0,c_m^0)\|_{L^1(\mathbb{R}^N)}\le C.
\end{array}
\right.
\end{equation}

An interesting difference with the purely mechanical model is that
when nutrients  are not enough the cells might die (or become quiescent), which is
represented by $\Phi (p,c)<0$ for $c$ small enough with $p$ given.
These cells, which are in principle in the center of the tumor, are
replaced by  cells from the boundary moving inwards by pressure
forces.  In particular the growth inequalities
\eqref{eq:rho_tpositive} cannot hold true here.

A typical choice is $\Psi(c,p) = c$, $\Phi(c,p) = \widetilde \Phi
(p) (c+\widetilde c_1)  - \widetilde c_2$ with $\widetilde c_i >0$
two constants and $\widetilde \Phi (p) $ a function having the same
properties as in the purely fluid mechanical model.

\subsection{Main results}


As in the case without nutrients, in the limit we obtain a free
boundary problem of Hele-Shaw type.

\begin{theorem}  Let $\Phi$, $\Psi$ satisfy \eqref{eq:assumptions.growth.nutrients},
and $\{(\varrho_m^0, c^0_m)\}$ satisfy the
hypotheses~\eqref{eq:assumptions.initial.data}
and~\eqref{eq:assumptions.initial.data.nutrients}. Then, after
extraction of subsequences,  the density $\varrho_m$, the nutrient
$c_m$ and  the pressure $p_m$ converge for all $T>0$ strongly in $L^1(Q_T)$  as
$m\to\infty$ to limits $\varrho_\infty,c_\infty,p_\infty \in BV(Q_T)$  that satisfy $0 \leq \varrho_\infty\leq 1$,  $0\leq
c_\infty\leq c_B$ and $0 \leq p_\infty \leq p_M$,  and
\begin{equation}
\label{eq:HSN} \left\{
\begin{array}{ll}
\displaystyle \partial_t\varrho_\infty  =\Delta p_\infty +\varrho_\infty
\Phi(p_\infty,c_\infty),\quad&\displaystyle\varrho_\infty(0)=\varrho^0,\\[4mm]
\displaystyle
\partial_t c_\infty =\Delta c_\infty - \varrho_\infty \Psi(p_\infty,c_\infty)\quad&c_\infty(0)=c^0,
\end{array}
\right.
\end{equation}
in a distributional sense, plus the relation $p_\infty\in
P_\infty(\varrho_\infty)$, where $P_\infty$ is the Hele-Shaw
graph~\eqref{eq:HS.graph}. \label{th:fbn}
\end{theorem}

When we say that $(\varrho_\infty,p_\infty,c_\infty)$ is a solution
to~\eqref{eq:HSN} in a distributional sense, we mean that
\begin{equation*}
\label{eq:distributional.sense}
\begin{array}{ll}
\displaystyle \iint_{Q} \big\{\varrho_\infty \psi_t+p_\infty\Delta
\psi+\varrho_\infty\Phi(p_\infty,c_\infty)\psi\big\}=-\int_{\mathbb{R}^N}\varrho^0\psi(\cdot,0),\\[4mm]
\displaystyle \iint_{Q} \big\{c_\infty \zeta_t+c_\infty\Delta
\zeta-\varrho_\infty\Psi(p_\infty,c_\infty)\zeta\big\}=-\int_{\mathbb{R}^N}c^0\zeta(\cdot,0),
\end{array}
\end{equation*}
for all test functions $\psi,\zeta\in C^\infty_0(\overline Q)$.

As mentioned above, for the system we lose the monotonicity
properties \eqref{eq:rho_tpositive}, and we are not able to prove the expected continuity of
$\varrho_\infty$ and $c_B-c_\infty$ in $L^1(\mathbb{R}^N)$. However, we are still able to
prove that  $\nabla p_\infty $ in $L^2(Q_T)$.
Therefore, the equation on the cell density can also be written as a
transport equation,
\begin{equation}
\label{eq:transport.system}
\partial_t  \varrho_\infty - {\rm div} \left(\varrho_\infty
\nabla p_\infty\right) = \varrho_\infty \Phi(p_\infty,c_\infty).
\end{equation}

If the initial data $\varrho_m^0$ are compactly supported uniformly
in $m$, we can control the supports of $\varrho_m(t)$ uniformly in
$m$ for each $t>0$. Hence, in the limit problem tumors propagate with
a finite speed (this is not true for the nutrients), and   a free
boundary shows up. This can be proved through a comparison argument
with the same barrier functions as in the scalar case. Moreover, the
constructed solution will fall within a class for which we prove
uniqueness.
\begin{theorem}
\label{thm:uniqueness.system} There is a unique triple $(\varrho, p,
c)$,  $\varrho,p,c_B-c\in L^\infty((0,\infty);L^1(\mathbb{R}^N)\cap
L^\infty(\mathbb{R}^N))$, $p\in P_\infty(\varrho)$,
satisfying~\eqref{eq:HSN} in the distributional sense and such that
for all $T>0$:
\begin{itemize}
\item[\rm (i)] $\rho(t)$ is uniformly compactly  supported for $t\in[0,T]$;
\item[(ii)] $|\nabla c|, |\nabla p| \in L^2 (Q_T)$;
\item[(iii)] $\partial_t p \in M^1 (Q_T)$, $\partial_t \varrho \in M^1(Q_T)$.
\end{itemize}

\end{theorem}

As in the purely mechanical model, we can write an equation for the
pressure $p_m$,
\begin{equation}
\label{eq:pressure_n}
\partial_t p_m=(m-1) p_m \Delta p_m +|\nabla p_m|^2 +(m-1)p_m\Phi(p_m,c_m).
\end{equation}
This suggests that  the complementarity formula
\begin{equation}
\label{eq:complementarity.system} p_\infty \left( \Delta p_\infty+
\Phi(p_\infty, c_\infty)\right)=0
\end{equation}
holds.  However, we have not been able to establish it rigorously. The reason is that, in contrast with the case without nutrients, we have failed to control
$\Delta p_m+ \Phi(p_m,c_m)$ from below  by means of a comparison argument, or to prove the strong convergence of $\nabla p_m$ in $L^2(Q_T)$ for all $T>0$. We leave
open the question to give conditions on  $\Phi(p,c) $ allowing  to
prove formula~\eqref{eq:complementarity.system}.

The rest of this section is devoted to prove the first of these theorems. The uniqueness result is
postponed to a later section.

\subsection{Estimates}
\noindent\textsc{$L^\infty$ bounds for $\varrho_m$, $p_m$, $c_m$. } The assumptions on the growth
functions~\eqref{eq:assumptions.growth.nutrients} imply, using
standard comparison arguments,
\begin{equation*} \label{eq:L-infty.bounds.nutrients}
0 \leq \varrho_m\leq \left(\frac{m-1}{m}p_M\right)^{1/(m-1)},\qquad
0\le p_m(x,t) \leq p_M,\qquad
0 < c_m < c_B.
\end{equation*}

\medskip

\noindent\textsc{$L^1$ bounds for $\varrho_m$, $p_m$, $c_m$. } Let $(\varrho,c)$  and $(\hat\varrho,\hat c)$ be two solutions of
the system for a fixed $m$, with corresponding pressures $p$ and $\hat p$. We
subtract the equation for $\hat \varrho$ from the equation for
$\varrho$, and multiply the resulting equation by $\mathop{\rm
sign}(\varrho-\hat\varrho)$, and do an analogous manipulation for the $c$
variable. We integrate by parts,  and obtain
$$
\left\{
\begin{array}{l}
\displaystyle \frac{d
}{dt}\|(\varrho-\hat\varrho)(t)\|_{L^1(\mathbb{R}^N)}\le
\int_{\mathbb{R}^N}\Big[\big(\varrho\Phi(p,c)-\hat \varrho\Phi(\hat p,\hat
c)\big)\mathop{\rm sign}(p-\hat p)\Big](t),
\\[10pt]
\displaystyle\frac{d }{dt}\|(c-\hat c)(t)\|_{L^1(\mathbb{R}^N)}\le-
\int_{\mathbb{R}^N}\Big[\big(\varrho \Psi(p,c)-\hat \varrho\Psi(\hat p,\hat
c)\big)\mathop{\rm sign}(c-\hat c)\Big](t).
\end{array}
\right.
$$
Let $\mu$ be a positive constant to be chosen later. We have
$$
\begin{array}{l}
\displaystyle\frac{d }{dt}\Big(\|(\varrho-\hat\varrho)(t)\|_{L^1(\mathbb{R}^N)} +  \mu
\|(c-\hat c)(t)\|_{L^1(\mathbb{R}^N)}\Big)\le\\[4mm]
\quad\displaystyle \underbrace{\int_{\mathbb{R}^N}\Big[\big(\varrho-\hat\varrho\big)\big(\Phi(p,c)\mathop{\rm sign}(\varrho-\hat\varrho)-\mu\Psi(p,c)\mathop{\rm sign}(c-\hat c)\big)\Big](t)}_{\mathcal{I}}\\[4mm]
\quad+\underbrace{\int_{\mathbb{R}^N}\Big[\hat\varrho\big((\Phi(p,c)-\Phi(p,\hat c))\mathop{\rm sign}(\varrho-\hat\varrho)
-\mu(\Psi(p,c)-\Psi(p,\hat c))\mathop{\rm sign}(c-\hat c)\big)\Big](t)}_{\mathcal{J}}\\[4mm]
\quad+\underbrace{
\int_{\mathbb{R}^N}
\Big[
\hat\varrho
\big(
(\Phi(p,\hat c)-\Phi(\hat p,\hat c))
\mathop{\rm sign}(\varrho-\hat\varrho)
-\mu(\Psi(p,\hat c)-\Psi(\hat p,\hat c))
\mathop{\rm sign}(c-\hat c)
\big)
\Big](t)
}_{\mathcal{K}}.
\end{array}
$$
Let $A=\{0\le p\le p_M,\, 0\le c\le c_B\}$, $\alpha=\min_{A}|\partial_p\Phi|>0$, $\beta=\max_{A}|\partial_p\Psi|$. Since $\partial_p \Phi<0$, if $\displaystyle 0<\mu\le \alpha/\beta$, there are constants $C$ independent of $m$ such that
$$
\begin{array}{l}
\mathcal{I}\le  \displaystyle\int_{\mathbb{R}^N}\Big[|\varrho-\hat\varrho|  \left(\Phi(p,c) + \mu \Psi(p ,c )\right)\Big](t)\leq
C\|(\varrho-\hat\varrho)(t)\|_{L^1(\mathbb{R}^N)},
\\[4mm]
\displaystyle\mathcal{J} \le  \| \hat \rho\|_{L^\infty(Q)} \left(\|
\partial_c\Phi\|_{L^\infty(A)}  + \mu \| \partial_c\Psi\|_{L^\infty(A)} \right)
\|(c - \hat c)(t) \|_{L^1(\mathbb{R}^N)} \leq C  \|(c -
\hat c)(t)\|_{L^1(\mathbb{R}^N)},\\[4mm]
\displaystyle\mathcal{K}\le\int_{\mathbb{R}^N} \Big[ \hat \varrho \big( -
\big|  \Phi(p,\hat c)- \Phi(\hat p,\hat c) \big| + \mu \big|
\Psi(p,\hat c)- \Psi(\hat p,\hat c) \big| \big)\Big](t)\le0.
\end{array}
$$
 We conclude that
$$
\frac{d }{dt}\Big(\|(\varrho-\hat\varrho)(t)\|_{L^1(\mathbb{R}^N)} +  \mu
\|(c-\hat c)(t)\|_{L^1(\mathbb{R}^N)}\Big) \le C \Big(\|(\varrho-\hat\varrho)(t)\|_{L^1(\mathbb{R}^N)}+  \mu
\|(c-\hat c)(t)\|_{L^1(\mathbb{R}^N)}\Big).
$$
Therefore,
$$
\|(\varrho-\hat\varrho)(t)\|_{L^1(\mathbb{R}^N)} + \mu \|(c-\hat c)(t)\|_{L^1(\mathbb{R}^N)} \le
\textrm{e}^{Ct}\Big(\|\varrho^0-\hat\varrho^0 \|_{L^1(\mathbb{R}^N)}+ \mu \|c^0-\hat
c^0\|_{L^1(\mathbb{R}^N)}\Big).
$$
This gives uniqueness, and choosing $(\varrho,c)=(\varrho_m,c_m)$, and
$(\hat\varrho, \hat c)= (0,c_B)$, we find the  uniform estimates
\begin{equation}\label{eqc:L1}
\|\varrho_m(t) \|_{L^1(\mathbb{R}^N)},\quad
\|c_m(t)-c_B\|_{L^1(\mathbb{R}^N)}\le K\textrm{e}^{Ct}.
\end{equation}

As in the purely fluid mechanical model, from the $L^\infty $ bounds
and \eqref{eqc:L1}, we conclude that, for $m>2$,
\begin{equation*}\label{eqc:pL1}
\| p_m(t)\|_{L^1(\mathbb{R}^N)} \le K\textrm{e}^{Ct}.
\end{equation*}

\medskip

\noindent\textsc{$L^1$ bounds on the derivatives of $\varrho_m$ and
$c_m$. } We differentiate the two equations of the system with
respect to time, multiply the first one by
$\text{sign}(\partial_t\varrho_m)$ and the second by
$\text{sign}(\partial_t c_m)$ and use Kato's inequality to obtain
\begin{equation}
\label{eq:system.derivatives}
\left\{
\begin{array}{l}
\displaystyle
\partial_t |\partial_t\varrho_m|-\Delta(m\varrho_m^{m-1}|\partial_t\varrho_m|)
\le |\partial_t\varrho_m|\Phi+\varrho_m\partial_p\Phi|\partial_t p_m|
+\varrho_m\partial_c\Phi \partial_t c_m\mathop{\rm sign}(\partial_t\varrho_m),\\[10pt]
\displaystyle
\partial_t |\partial_t c_m|-\Delta(|\partial_t c_m|)\le -\partial_t\varrho_m
\Psi\mathop{\rm sign}(\partial_t c_m)-\varrho_m\partial_p\Psi \partial_t
p_m\mathop{\rm sign}(\partial_t c_m)-\varrho_m\partial_c\Psi|\partial_t
c_m|.
\end{array}
\right.
\end{equation}
We integrate in space and  add the two equations to obtain, using
the monotonicity properties of the growth functions,
$$
\begin{gathered}
\frac{d}{dt}\Big(\|\partial_t\varrho_m(t)\|_{L^1(\mathbb{R}^N)}+\mu \|\partial_t
c_m(t)\|_{L^1(\mathbb{R}^N)}\Big) \le\\[4mm]
 \int_{\mathbb{R}^N} \left( |\partial_t \varrho_m(t)| \big(\Phi +
\mu \Psi\big)(0,c_B) + \| \varrho_m \|_{L^\infty(Q)}
\|\partial_c\Phi\|_{L^\infty(A)}|\partial_t c_m(t)| \right),
\end{gathered}
$$
and thus, thanks to assumption \eqref{eq:assumptions.initial.data.nutrients},
\begin{equation*}\label{eqc:tderivatives}
\|\partial_t\varrho_m(t)\|_{L^1(\mathbb{R}^N)}+\|\partial_t c_m(t)\|_{L^1(\mathbb{R}^N)}\le
C\textrm{e}^{C t}
\Big(\|(\partial_t\varrho)^0\|_{L^1(\mathbb{R}^N)}+|(\partial_t c)^0\|_{L^1(\mathbb{R}^N)}\Big) =
K\textrm{e}^{Ct}.
\end{equation*}

We can estimate  the space derivatives in the same way,  and arrive to
\begin{equation*}\label{eqc:xderivatives}
\|\partial_{x_i}\varrho_m(t)\|_{L^1(\mathbb{R}^N)}+\|\partial_{x_i}
c_m(t)\|_{L^1(\mathbb{R}^N)} \le K\textrm{e}^{Ct}, \qquad i=1,\dots,N.
\end{equation*}

\medskip

\noindent\textsc{Estimates on the derivatives of $p_m$. } From the
first equation in \eqref{eq:system.derivatives},  using the strict
monotonicity of $\Phi$ with respect to $p$ and the $BV$ estimates
that we have already proved, we get
$$
\int_{\mathbb{R}^N}|\partial_t\varrho_m (T)|+\alpha
\iint_{Q_T}\varrho_m|\partial_t p_m|\le K\textrm{e}^{CT},
$$
where, as above, $\alpha=\min_{A}|\partial_p\Phi|>0$.
Therefore,
\begin{equation*}\label{eqc:tderivativesp}
\|\partial_t p_m\|_{L^1(Q_T)}\le\iint_{Q_T\cap\{\varrho\le1/2\}}m\varrho_m^{m-2}|\partial_t\varrho_m|+
2\iint_{Q_T\cap\{\varrho_m\ge1/2\}}\varrho_m|\partial_t
p_m|\le C(T).
\end{equation*}

The control on the space derivatives follows by a similar argument,
and we obtain
\begin{equation*}\label{eqc:xderivativesp}
\|\partial_{x_i} p_m\|_{L^1(Q_T)} \le C(T),\quad i=1,\dots,N.
\end{equation*}

\medskip

\noindent\textsc{Convergence and identification of the limit. } The
above estimates give strong convergence in $L^1_{\rm loc}(Q)$ for
$\varrho_m$, $c_B-c_m$ and $p_m$. Our uniform control of the $L^1$
norms of these quantities in $m$ and  $t\in[0,T]$ shows that the
mass in thin initial strips is uniformly in $m$ small. Therefore, to
get global convergence in $L^1(Q_T)$, we just need to control the
tails. Proceeding as in the case without nutrients, we get
$$
\frac{d}{dt}\int_{\mathbb{R}^N}\varrho_m(t)\varphi_R\le
CR^{-2}+\Phi(0,c_B)\int_{\mathbb{R}^N}\varrho_m(t)\varphi_R,
$$
from where the tail control for $\varrho_m$, and then for $p_m$, follows easily.

As for the nutrients, we have
$$
\frac{d}{dt}\int_{\mathbb{R}^N}(c_B-c_m(t))\varphi_R\le
R^{-2}\|\Delta\varphi\|_{L^\infty(\mathbb{R}^N)}\int_{\mathbb{R}^N}(c_B-c_m(t))
+\Psi(0,c_B)\int_{\mathbb{R}^N}\varrho_m(t)\varphi_R.
$$
But we already have a tail control for $\varrho_m$. Hence, for all $R$ and $m$ large enough,
$$
\int_{\mathbb{R}^N}(c_B-c_m(t))\varphi_R\le
\int_{\mathbb{R}^N}(c_B-c^0(t))
\varphi_R+\int_{\mathbb{R}^N}|c^0(t)-c^0_m(t)|\varphi_R+CR^{-2}+\Psi(0,c_B)\varepsilon,
$$
which immediately implies a tail control for $c_B-c_m$.

The strong convergence in $L^1(Q_T)$ is enough  to pass to
the limit and recover the system \eqref{eq:HSN} in a distributional sense with the limiting
pressure graph relation.

\medskip

\noindent\textsc{$L^2$ bounds for $\nabla p_m$, $\nabla c_m$. }
The bound for the gradients of the pressures follows from equation
\eqref{eq:pressure_n} written in the form
$$
\partial_t p_m =\frac{m-1}2\Delta p_m^2-(m-2)|\nabla p_m|^2+(m-1)p_m\Phi(p_m,c_m).
$$
Indeed, integrating in $Q_T$ we obtain
$$
\|\nabla p_m\|_{L^2(Q_T)}^2\le\frac{m-1}{m-2}\Phi(0,c_B) \|p_m\|_{L^1(Q_T)}
+\frac{1}{m-2}\|p^0_m\|_{L^1(\mathbb{R}^N)}.
$$
This is enough to prove that the limit satisfies the transport
equation~\eqref{eq:transport.system}.

The estimate for the gradients of the nutrients, which is needed to
prove that  the limit falls within the uniqueness class, is even
easier. We just have to multiply the equation for the nutrients by
$c_m$ and integrate by parts in $Q_T$. Since $\Psi\ge0$, we obtain
$$
\|\nabla c_m\|_{L^2(Q_T)}\le\|c^0_m\|_{L^2(\mathbb{R}^N)}/\sqrt{2}.
$$


\section{Uniqueness for the limit system with nutrients}
\label{sec:nutrient_unique}
\setcounter{equation}{0}

We meet two difficulties. On the one hand,  as in the case without
nutrients, $p$ is not a Lipschitz, single-valued function of
$\varrho$. On the other hand, the density of nutrients is not
compactly supported. Hence we can not restrict to the case of a
bounded domain. To take care of the lack of compact support of $c$,
we will use an idea from~\cite{Benilan-Crandall-Pierre-1984}.

\medskip

\noindent\textsc{Dual system. }
Let $(\varrho_1, p_1, c_1)$ and $(\varrho_2, p_2, c_2)$ be two
solutions with the same initial data and fix a final time $T$. Since
they satisfy \eqref{eq:HSN} in the sense of distributions, we have
\begin{equation}
\label{eq:uniqueness.sense.distributions}
\begin{cases}
\displaystyle\iint_{Q_T}\left\{
(\varrho_1-\varrho_2)   \partial_t \psi +  (p_1-p_2) \Delta \psi + (
\varrho_1 \Phi(p_1,c_1) - \varrho_2 \Phi(p_2,c_2)  )  \psi \right\}
=0 ,
\\[10pt]
\displaystyle\iint_{Q_T}\left\{ (c_1-c_2)
\partial_t \zeta    +  (c_1-c_2) \Delta \zeta - ( \varrho_1
\Psi(p_1,c_1) - \varrho_2 \Psi(p_2,c_2)  )\zeta \right\} = 0 ,
\end{cases}
\end{equation}
for any pair of test functions $\psi,\zeta\in\mathcal{D}(\overline Q)$ such
that  $\psi(\cdot,T) = \zeta(\cdot,T)=0$.  Denoting $\Phi_1=\Phi(p_1,c_1)$ and
$\Psi_1=\Psi(p_1,c_1)$, we can rewrite the above equations as
$$
\begin{cases}
\displaystyle\iint_{Q_T} \left\{(\varrho_1- \varrho_2+p_1-p_2)
\left(A\partial_t \psi  +  B \Delta \psi + A\Phi_1 \psi - CB \psi
\right) + (c_1-c_2) D \psi\right\} =0 ,
\\[10pt]
\displaystyle \iint_{Q_T}  \left\{(c_1-c_2)  \left(
\partial_t \zeta  +  \Delta \zeta - E  \zeta \right)  - (\varrho_1-
\varrho_2+p_1-p_2) \left(A\Psi_1   - FB \right) \zeta \right\}= 0,
\end{cases}
$$
where
$$
\begin{array}{ll}
0 \leq \displaystyle
A=\frac{\varrho_1-\varrho_2}{\varrho_1-\varrho_2 +p_1-p_2} \leq
1,\ & 0 \leq \displaystyle
B=\frac{p_1-p_2}{\varrho_1-\varrho_2+p_1-p_2} \leq  1,
\\[12pt]
0 \leq \displaystyle C= \;
-\varrho_2\frac{\Phi(p_1,c_1)-\Phi(p_2,c_1)}{p_1- p_2} \le \nu_1,  \
&\displaystyle 0\leq D= \varrho_2
\frac{\Phi(p_2,c_1)-\Phi(p_2,c_2)}{c_1- c_2}\le \nu_2,
\\[12pt]
\displaystyle0\leq E= \varrho_2
\frac{\Psi(p_2,c_1)-\Psi(p_2,c_2)}{c_1- c_2} \le \nu_3,\ &0 \leq
\displaystyle F= \;
-\varrho_2\frac{\Psi(p_1,c_1)-\Psi(p_2,c_1)}{p_1- p_2}\le \nu_4.
\end{array}
$$
Adding these two equations we get
\begin{equation}
\label{eq:added.weak.system} \displaystyle\iint_{Q_T}
\Big(\big(\varrho_1-  \varrho_2+p_1-p_2\big) \mathcal{A}(\psi,\zeta)
+ (c_1-c_2)\mathcal{B}(\psi,\zeta)\Big) =0 ,
\end{equation}
where
\begin{equation}
\label{eq:added.weak.system.coefficients} \left\{\begin{array}{l}
\displaystyle\mathcal{A}(\psi,\zeta)=A\partial_t \psi +  B \Delta
\psi
+  \left(A\Phi_1 - CB\right) \psi -\left(A\Psi_1  -FB\right)  \zeta,\\[10pt]
\displaystyle\mathcal{B}(\psi,\zeta)=\partial_t \zeta   +  \Delta
\zeta - E\zeta  +D \psi.
\end{array}
\right.
\end{equation}

Let $G$, $H$ be any non-negative functions in $\mathcal{D}'(Q_T)$.
If the \emph{dual system}
$$
\mathcal{A}(\psi,\zeta)=AG,\qquad \mathcal{B}(\psi,\zeta)=H, \quad
\psi(\cdot,T)=0,\quad \zeta(\cdot,T)=0,
$$
admits a smooth solution $(\psi,\zeta)$ with a suitable decay, we
may  use $\psi$ and $\zeta$ as test functions to obtain
$$
\displaystyle\iint_{Q_T} \big((\varrho_1- \varrho_2)G +
(c_1-c_2)H\big)=0 ,
$$
from where uniqueness would follow. Unfortunately, the coefficients
in the equations defining the dual system are not smooth. Even
worse, $A$ an $B$ may vanish: the system is not uniformly parabolic,
and a delicate regularization procedure is required to fulfill our
plan.

\medskip

\noindent\textsc{Regularized dual system. } Given  $G$, $H$ as
above, let $R_0>0$ be such that the supports of $G(\cdot,t)$,
$H(\cdot,t)$, $\varrho_1(t)$, $\varrho_2(t)$, $t\in(0,T)$ are
contained in $B_{R_0}(0)$. Given any $R>R_0+1$, we introduce the
dual system with regularized  coefficients, posed in
$Q_{R,T}=B_R(0)\times(0,T)$,
\begin{equation}
\label{eqn:system_dual}
\left\{
\begin{array}{l}
\displaystyle
 \partial_t \psi_{n,R} +  \frac{B_n}{A_n} \Delta \psi_{n,R}
+  \left(\Phi_{1,n} - \frac{C_nB_n}{A_n}\right) \psi_{n,R} -  \left(
\Psi_{1,n}  - \frac{F_nB_n}{A_n}\right) \zeta_{n,R}= G,
 \\[10pt]
\displaystyle
\partial_t \zeta_{n,R}    +  \Delta \zeta_{n,R} - E_n  \zeta_{n,R} +D_n  \psi_{n,R} = H,
\end{array}
\right.
\end{equation}
with final and boundary data given by
$$
\psi_{n,R}(\cdot,T)=\zeta_{n,R}(\cdot,T)=0,\qquad
\psi_{n,R}=\zeta_{n,R}=0 \text{ in }\partial B_{R}(0)\times(0,T),
$$
where, thanks to the hypotheses on the solutions, the $L^2(Q_{T})$
error in the approximation of each coefficient is $O(1/n)$, and
\begin{equation*}
\left\{\begin{array}{l}
1/n<A_n, B_n\le 1,\quad \ 0 \leq C_n,D_n,E_n,F_n <M,\quad |\Phi_{1,n}|,\ |\Psi_{1,n}|<M,\\[3mm]
\| \partial_t C_n \|_{L^1(Q_{T})},\ \| \partial_t F_n
\|_{L^1(Q_{T})},\ \| \nabla \Phi_{1,n}  \|_{L^2(Q_{T})},\ \| \nabla
\Psi_{1,n}  \|_{L^2(Q_{T})} \leq K.
\end{array}
\right.
\end{equation*}
Our aim is to  use $\psi_{n,R}$ and $\zeta_{n,R}$, suitably extended
by zero, as test functions. Since these functions have non-zero
derivatives at the boundary,  we find yet another difficulty: the
lack of smoothness at the boundary of the ball $B_R(0)$, which
requires an extra regularization.

\medskip

\noindent\textsc{Avoiding the lack of smoothness at the boundary of the ball. }  For
$0<\varepsilon<1/2$, we consider a family of cut-off functions
$\eta_\varepsilon\in C^\infty_0(\mathbb{R}^N)$ such that $0\le
\eta_\varepsilon\le 1$,
$$
\begin{array}{l}
\eta_\varepsilon(x)=1 \text{ if } |x|<R-2\varepsilon,\qquad \eta_\varepsilon=0
\text{ if } |x|>R-\varepsilon,\\[10pt]
\|\nabla\eta_\varepsilon\|_\infty \le \mathcal{C}/\varepsilon,\qquad
\|\Delta\eta_\varepsilon\|_\infty \leq \mathcal{C}/\varepsilon^2.
\end{array}
$$

Now put $\psi=\eta_\varepsilon\psi_{n,R}$,
$\zeta=\eta_\varepsilon\zeta_{n,R}$  in the weak
formulation~\eqref{eq:added.weak.system}--\eqref{eq:added.weak.system.coefficients}
to get
\begin{equation*}
\label{eq:terms.approximation.system} \displaystyle\iint_{Q_T}
\eta_\varepsilon\big((\varrho_1- \varrho_2)G + (c_1-c_2)H\big) =-
\mathcal{J}_{n,R,\varepsilon}-\mathcal{K}_{n,R,\varepsilon}-\sum_{i=1}^5I^i_{n,R,\varepsilon},
\end{equation*}
where
$$
\begin{array}{l}
\mathcal{J}_{n,R,\varepsilon}=\displaystyle\iint_{Q_{R,T}}\left\{(p_1-p_2)
\big(2\nabla\eta_\varepsilon\cdot\nabla\psi_{n,R}
+\psi_{n,R}\Delta\eta_\varepsilon\big)\right\},\\[10pt]
\mathcal{K}_{n,R,\varepsilon}=\displaystyle\iint_{Q_{R,T}}\left\{(c_1-c_2)
\big(2\nabla\eta_\varepsilon\cdot\nabla\zeta_{n,R}
+\zeta_{n,R}\Delta\eta_\varepsilon\big)\right\},\\[10pt]
I^1_{n,R,\varepsilon}=\displaystyle\iint_{Q_{R,T}}\left\{(\varrho_1-\varrho_2+
p_1-p_2)\eta_\varepsilon B_n\left(1-\frac{A}{A_n}\right)
\left(\Delta\psi_{n,R}-C_n\psi_{n,R}+F_n\zeta_{n,R}\right)\right\},\\[10pt]
I^2_{n,R,\varepsilon}=\displaystyle\iint_{Q_{R,T}}\left\{(\varrho_1-\varrho_2+p_1-p_2)
\eta_\varepsilon (B-B_n)\left(\Delta\psi_{n,R}-C_n\psi_{n,R}+F_n\zeta_{n,R}\right)\right\},\\[10pt]
I^3_{n,R,\varepsilon}=-\displaystyle\iint_{Q_{R,T}}\big\{(\varrho_1-\varrho_2+p_1-p_2)
\eta_\varepsilon B [\psi_{n,R}(C-C_n)+\zeta_{n,R}(F-F_n) ] \big\},\\[10pt]
I^4_{n,R,\varepsilon}=\displaystyle\iint_{Q_{R,T}}\big\{(c_1-c_2)
\eta_\varepsilon\zeta_{n,R} [(E-E_n)+(D-D_n)] \big\},\\[10pt]
I^5_{n,R,\varepsilon}=\displaystyle\iint_{Q_{R,T}}\left\{(\varrho_1-\varrho_2+p_1-p_2)\eta_\varepsilon
A\left((\Phi_1-\Phi_{1,n})\psi_{n,R}-
(\Psi_1-\Psi_{1,n})\zeta_{n,R}\right)\right\}.
\end{array}
$$
The term $\mathcal{J}_{n,R,\varepsilon}$ vanishes, because the
supports of $p_1(t)$, $p_2(t)$, $t\in(0,T)$ are contained in
$B_{R_0}(0)$. Hence we can drop it.

We will now take the limit $\varepsilon\to0$. The only difficult
term at  this step is
$\mathcal{K}_{n,R,\varepsilon}$.
We write
$$
| \mathcal{K}_{n,R,\varepsilon} | \leq \mathcal{C}\displaystyle\int_0^T
\int_{R-2\varepsilon<|x|<R}  | c_1-c_2 |
\left(\frac{|\nabla\zeta_{n,R}|
}{\varepsilon}+\frac{|\zeta_{n,R}|}{\varepsilon^2}\right).
$$
Since $\zeta_{n,R}=0$ on $\partial B_R(0)$, we have
$$
\sup_N|\zeta_{n,R}|\le 2\varepsilon |\nabla\zeta_{n,R}|, \qquad
N=\{R-2\varepsilon<|x|<R,\ 0<t<T\},
$$
and, if $\nu$ is the unit outwards normal to $\partial B_{R}(0)$,  then
$$
\lim_{\varepsilon\downarrow0}\left(\sup_N|\nabla\zeta_{n,R}|\right)=
\sup_{|x|=R,0<t<T}|\nabla\zeta_{n,R}|=
\sup_{|x|=R,0<t<T}|\partial_\nu\zeta_{n,R}|.
$$
Therefore, since $c_1$ and $c_2$ are bounded,
$$
\limsup_{\varepsilon\downarrow0} | \mathcal{K}_{n,R,\varepsilon} |
\le \mathcal{C}R^{N-1}\sup_{|x|=R,
0<t<T}|\partial_\nu\zeta_{n,R}|.
$$

\medskip

\noindent\textsc{Limit $n\to\infty$. Uniform bounds. }  Proving that $\lim_{n\to\infty}I^i=0$, $i=1,\dots,5$, is easy, thanks to the
following estimates.

\begin{lemma}  There are constants $\mathcal{C}_1$, $\mathcal{C}_2$, depending on $T$, $G$ and
$H$, but not on $n$ and $R$, such that
\begin{equation*}
\label{D:max}
\begin{array}{l}
\| \psi_{n,R}\|_{L^\infty(Q_{R,T})},\ \| \zeta_{n,R}\|_{L^\infty(Q_{R,T})}  \leq \mathcal{C}_1,\\[4mm]
\left\| \left(B_n/A_n\right)^{1/2} \left( \Delta \psi_{n,R} -C_n
\psi_{n,R}+F_n\zeta_{n,R}\right)\right\|_{L^2(Q_{R,T})} \leq
\mathcal{C}_2.
\end{array}
\end{equation*}
\end{lemma}

\noindent\emph{Proof. }  \emph{$L^\infty$ bounds. } An easy
computation shows that the functions
$$
M(t) = \max_{|x| \leq R} | \psi_{n,R}(x,t) | ,  \qquad N(t) =
\max_{|x| \leq R}  | \zeta_{n,R}(x,t) |,
$$
satisfy the differential inequalities
$$
\left\{
\begin{array}{l}
- M'(t) + \alpha_n(t) M(t)  \leq \mathcal{C} (1+ M(t) + N(t)) + q_n(t)
\alpha_n(t) N(t),
 \\[10pt]
- N'(t)  \leq  \mathcal{C}(1+ M(t) ),
 \\[10pt]
 M(T)=N(T)=0,
\end{array}
\right.
$$
with $\alpha_n(t)=  \frac{C_nB_n}{A_n}(x_n,t) \geq 0$,
$q_n(t)=\frac{F_n}{C_n}(x_n,t)$,  where $x_n$ is  the point where
the maximum defining $M(t)$ is achieved. Thanks to the assumptions
on $\partial_p \Phi$ and $\partial_p\Psi$, there is a constant
$\overline q>0$ such that $0\le q_n(t)\le\overline q$.

Let now
$$
Q(t) = \max\big(  M(t), \overline q N(t) \big).
$$
As a combination of the equations on $M$ and $N$, we find that
\begin{equation}
\label{eq:differential.inequality.Linfty.estimate}
 - Q'(t)  \leq
\mathcal{C} (1+ Q(t)), \qquad Q(T)=0.
\end{equation}
Indeed, assume that $\max\big(  M(t_0), \overline q N(t_0)
\big)=M(t_0)$. Then, $q_n(t_0) N(t_0) \leq M(t_0)$ and the bad
$\alpha_n$ terms cancel. If, on the contrary, $\max\big(  M(t_0),
\overline q N(t_0) \big)=N(t_0)$, then the result follows by the
harmless equation on $N(t)$.

From~\eqref{eq:differential.inequality.Linfty.estimate} we have
$Q(t)\le 1-\textrm{e}^{\mathcal{C}{T-t}}$, hence the result.

\medskip

\noindent\emph{Estimates on the Laplacian. } We multiply the second equation
in~\eqref{eqn:system_dual} by $\partial_t\zeta_{n,R}$ and integrate
in space and time to obtain
$$
\begin{array}{rcl}
\displaystyle\|\partial_t\zeta_{n,R}\|_{L^2(Q_{R,T})}^2+\frac12\|\nabla
\zeta_{n,R}(t)\|_{L^2(B_R(0))}^2&=&\displaystyle\iint_{Q_{R,T}}(E_n\zeta_{n,R}-D_n\psi_{n,R}+H)\partial_t\zeta_{n,R}\\[4mm]
&\le&\displaystyle C\|\partial_t\zeta_{n,R}\|_{L^2(Q_{R,T})}.
\end{array}
$$
Hence $\|\partial_t\zeta_{n,R}\|_{L^2(Q_{R,T})}$ is uniformly
bounded.

We now multiply the first equation in \eqref{eqn:system_dual} by
$\Delta \psi_{n,R} -C_n \psi_{n,R}+F_n\zeta_{n,R}$ and integrate in
space and time. We get
$$
\begin{array}{l}
\displaystyle\frac12\|\nabla
\psi_{n,R}(t)\|_{L^2(B_R(0))}^2+\int_t^T\int_{B_R(0)}\frac{B_n}{A_n}
\Big(\Delta
\psi_{n,R} -C_n \psi_{n,R}+F_n\zeta_{n,R}\Big)^2\\[4mm]
\
\displaystyle=-\int_{B_R(0)}\left(\frac{C_n\psi_{n,R}^2}2+F_n\psi_{n,R}\zeta_{n,R}\right)(t)\\[4mm]
\quad\displaystyle+\int_t^T\int_{B_R(0)}
\left(-\frac{\psi_{n,R}^2\partial_t C_n}2+
\psi_{n,R}\zeta_{n,R}\partial_t F_n+F_n\psi_{n,R}\partial_t\zeta_{n,R}\right)\\[4mm]
\quad \displaystyle
+\int_t^T\int_{B_R(0)}\Big(-\Phi_{1,n}|\nabla\psi_{n,R}|^2+\psi_{n,R}\nabla\Phi_{1,n}\cdot
\nabla\psi_{n,R}+ C_n\Phi_{1,n}\psi_{n,R}^2-
\Phi_{1,n}F_n\psi_{n,R}\zeta_{n,R}\Big)\\[4mm]
\quad \displaystyle +\int_t^T\int_{B_R(0)}\Big(
-\Psi_{1,n}\nabla\zeta_{n,R}\cdot\nabla\psi_{n,R}-\zeta_{n,R}\nabla\Psi_{1,n}\cdot
\nabla\psi_{n,R}-C_n\Psi_{1,n}\psi_{n,R}\zeta_{n,R}+F_n\psi_{n,R}\zeta_{n,R}^2 \Big)\\[4mm]
\quad \displaystyle +\int_t^T\int_{B_R(0)}\Big(\psi_{n,R}\Delta
G-C_nG\psi_{n,R}+F_nG\zeta_{n,R} \Big)\\[4mm]
\ \displaystyle\le
K\left(1-t+\int_t^T\|\nabla\psi_{n,R}(s)\|_{L^2(B_R(0))}\,ds\right).
\end{array}
$$
Using the Gr\"{o}nwall lemma we get $\sup_{t\in(0,T)}\|\nabla
\psi_{n,R}(t)\|_{L^2(B_R(0))}\le C$, and then the desired estimate
\qed

In order to pass to the limit in the $\mathcal{K}$ term we need a good enough estimate,  uniform in $n$, for
the decay of the normal derivative.

\begin{lemma}
As $R\to\infty$,
\begin{equation}
\label{eq:estimate.normal.derivative} \sup_{|x|=R,
0<t<T}|\partial_\nu\zeta_{n,R}|= o(R^{1-N})
\end{equation}
uniformly in $n$.
\end{lemma}

\noindent\emph{Proof. } As a first step we obtain a more precise
estimate for the size of $\zeta_{n,R}$. We first notice that, since
$E_n,D_n\ge0$, $D_n$ and $H$ are bounded and compactly supported in
$B_{R_0}(0)$, and $\psi_{n,R}$ is bounded, then
$$
\partial_t \zeta_{n,R}    +  \Delta \zeta_{n,R}=E_n  \zeta_{n,R} -D_n  \psi_{n,R} + H\ge
-\mathcal{C} \mathds{1}_{B_{R_0}(0)}
$$
for some $\mathcal{C}>0$. On the other hand, given a fixed $\tau>0$,
if $K>0$ is large enough, the function
$$
Z(x,t)=K\textrm{e}^{-\lambda
t}\frac{\textrm{e}^{-|x|^2/(4(T-t+\tau))}}{(4\pi (T-t+\tau))^{N/2}}
$$
satisfies
$$
\begin{array}{rcl}
\displaystyle Z_t+\Delta Z&=&\displaystyle-\lambda
K\textrm{e}^{-\lambda
t}\frac{\textrm{e}^{-|x|^2/(4(T-t+\tau))}}{(4\pi
(T-t+\tau))^{N/2}}\le-\mathds{1}_{B_{R_0}(0)}\lambda K
\frac{\textrm{e}^{-R_0^2/(4\tau)}}{(4\pi (T+\tau))^{N/2}}
\\[4mm]
&\le&\displaystyle-\mathcal{C} \mathds{1}_{B_{R_0}(0)}.
\end{array}
$$
Moreover, $Z(x,T)=0$, and $Z(x,t)\ge0$ for $|x|=R$ and $0<t<T$.
Therefore, comparison yields
$$
0\le \zeta_{n,R}\le Z \qquad\text{in } Q_{R,T}.
$$

If $N\ge 3$, we consider a function $g=g(x)$ defined by
$$
g(x)=\frac{d}{|x|^{N-2}}+e,
$$
where $e$, $d$ satisfy
$$
\frac{d}{(R-1)^{N-2}}+e=K\frac{\textrm{e}^{-R^2/(4(T+\tau))}}{(4\pi
\tau)^{N/2}},\qquad \frac{d}{R^{N-2}}+e=0.
$$
Note that $\Delta g=0$ on $\{R-1<|x|<R\}$. Moreover,
$$
\begin{array}{ll}
g(x)\ge \zeta_{n,R}(x,t),\qquad&|x|=R-1,\ 0<t<T,\\[4mm]
g(x)=\zeta_{n,R}(x,t)=0,\qquad&|x|=R,\ 0<t<T,\\[4mm]
g(x)\ge\zeta_{n,R}(x,T)=0,\qquad&R-1<|x|<R.
\end{array}
$$
Therefore, $g(x)\ge\zeta_{n,R}(x,t)$ for all $R-1<|x|<R$, $0<t<T$.
Since $g(x)=0=\zeta_{n,R}(x,t)$ for $|x|=R$,  $0<t<T$, we conclude
that
$$
\partial_\nu(g-\zeta_{n,R})(x,t)\quad\text{for } |x|=R,\ 0<t<T,
$$
and hence
$$
\sup_{|x|=R, 0<t<T}|\partial_\nu\zeta_{n,R}|\le \sup_{|x|=R, 0<t<T}|\partial_\nu g|.
$$
The result then follows from the following estimate,
$$
\begin{array}{rcl}
\displaystyle R^{N-1}|\partial_\nu
g|(x)&=&\displaystyle(N-2)d=\frac{(N-2)K\textrm{e}^{-R^2/(4(T+\tau))}}{(4\pi
\tau)^{N/2}}
\left(\frac{1}{(R-1)^{N-2}}-\frac{1}{R^{N-2}}\right)^{-1}\\[4mm]
&=&\displaystyle o(1)\quad\text{as }R\to\infty.
\end{array}
$$

An analogous computation with $g(x)=d\log|x|+e$ if $N=2$, and $g(x)=d|x|+e$ if $N=1$, lead to the same estimate~\eqref{eq:estimate.normal.derivative}.
\qed

\medskip

\noindent\textsc{Limit $R\to\infty$ and conclusion. } Passing to the
limit in $R$, we finally obtain that
$$
\iint_{Q_T} \big( (\varrho_1- \varrho_2)G+(c_1-c_2)H \big)=0.
$$
Hence  $\varrho_1=\varrho_2$ and $c_1=c_2$. Uniqueness for $p$ then
follows just taking  $\psi=p_1-p_2$ as test function
in~\eqref{eq:uniqueness.sense.distributions}.

\section*{Appendix: Examples for the purely mechanical model}
\renewcommand{\theequation}{A.\arabic{equation}}

\renewcommand{\thesection}{A}
\setcounter{equation}{0}
\setcounter{theorem}{0}

\subsection{Tumor spheroids}
\label{sec:spheroids}
\begin{figure}[h!]
\centerline{
 \includegraphics[width=12cm]{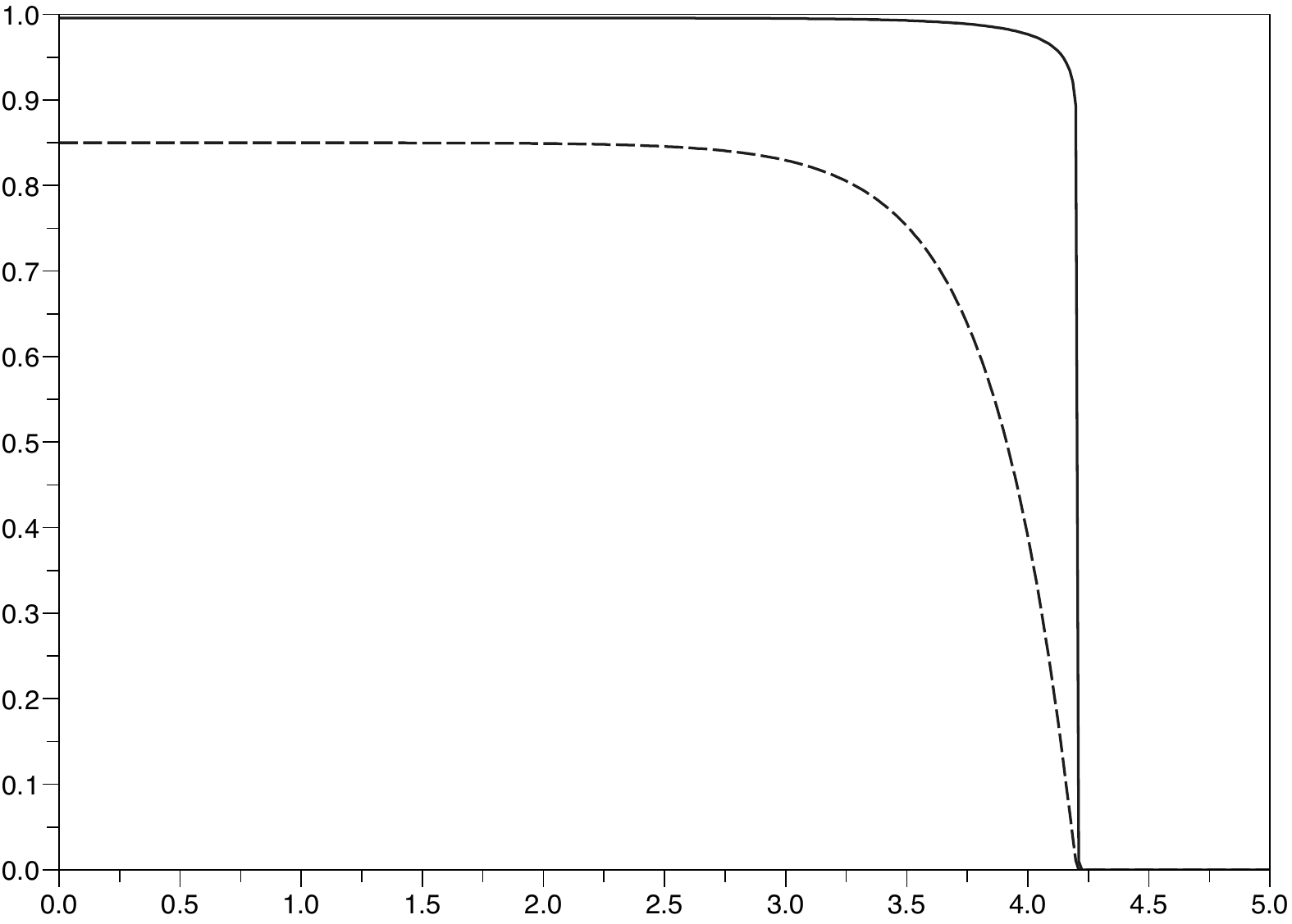}}
\vspace{-5mm} \caption{\emph{Traveling wave. } A traveling wave
solution to the mechanical model \eqref{model1},
\eqref{eq:def.pressure} in one dimension with $m=40$. The upper
continuous line is $\varrho$; the bottom dashed line is $p$. Here
$p_M=.85$.} \label{fig:tw}
\end{figure}

A typical application of the Hele-Shaw equations is to describe
tumor spheroids \cite{Bru,
byrne-chaplain,byrne-drasdo,cui_escher,friedman,friedman_hu,Lowengrub_survey}.
When nutrients are ignored, the tumor is assumed to fill a ball
centered at $0$,
$$
\Omega(t) := \{ p_\infty(t) >0\}  =\{ \varrho_\infty(t) =1 \} =
B_{R(t)}(0).
$$
The radius $R(t)$ of this ball is computed according to the
geometric motion rules~\eqref{eq:elliptic.equation}
and~\eqref{eq:geometric.motion.law};  that is,  we consider the
unique (and thus radially symmetric) solution to
\begin{equation} \label{apts_p}
- \Delta p_\infty(t) = \Phi(p_\infty(t)) \quad \text{in }
B_{R(t)}(0), \quad \qquad p_\infty(R(t), t )=0,
\end{equation}
and evolve the radius according to
\begin{equation} \label{apts_R}
R'(t)=V= |\nabla p_\infty(R(t), t)|.
\end{equation}
Then, we consider  $ \varrho_\infty $  defined as
\begin{equation} \label{apts_rho}
\varrho_\infty (t)= \mathds{1}_{B_{R(t)}(0)} .
\end{equation}
This is indeed a correct solution to our model.

\begin{theorem}Let $R(0)=R^0$ be given.
Problem \eqref{apts_p}--\eqref{apts_rho} defines a unique dynamic
$R(t)$, $\varrho_\infty (t)$, $p_\infty (t)$, which turns out to be
the unique solution to the Hele-Shaw limit problem~\eqref{eq:HS}
with initial data $\varrho_\infty^0=\mathds{1}_{B_{R^0}(0)}$. For
long times it approaches a \lq traveling wave' solution with a
limiting speed independent of the dimension,
\begin{equation} \label{apts_speed}
 R'(t) \underset{ t \to \infty }{\longrightarrow} \sqrt{2 Q(p_M)},
\qquad Q(p)=\int_0^p \Phi(q)dq.
\end{equation}
\label{th:ts}
The limit profile can also be calculated and is one-dimensional.
\end{theorem}

For several more elaborate one dimensional models, it is also
possible  to compute the traveling waves which
define the asymptotic shape for large times; see \cite{TVCVDP}.

\noindent \emph{Proof. } \emph{Existence. }  A unique solution to
\eqref{apts_p}  can be defined for every fixed $R(t)$, because the elliptic
problem \eqref{apts_p} comes with a decreasing nonlinearity. Indeed,
shooting or iterative methods apply and give a solution
with smooth dependency on the parameter $R(t)$. Therefore the
differential equation \eqref{apts_R} can be solved thanks to the
Cauchy-Lipschitz theorem.

\medskip

\noindent\emph{Equivalence with the Hele-Shaw problem. } To show
that this geometric  method gives the solution to  \eqref{eq:HS}, we
compute, in the distributional sense, that
$$
\partial_t \varrho_\infty = R' (t) \delta_{\partial \Omega(t)}.
$$
Also, because $p_\infty$ satisfies $\Delta p_\infty= \Phi(p_\infty)$
in $\Omega(t)$ and since $p_\infty$ vanishes out of $\Omega(t)$, we
have, still in the  distributional sense in ${\mathbb R}^N$, with
$\nu$ the unit outward normal to $ \Omega(t)$,
$$
\Delta p_\infty + \varrho_\infty \Phi(p_\infty) = - \delta_{\partial
\Omega(t)} \partial_\nu p_\infty =  - \delta_{\partial \Omega(t)}
|\nabla p_\infty(R(t), t)|,
$$
and the last inequality follows from the Dirichlet  boundary
condition. This means that both formulations, geometric motion of
the domain or  equation \eqref{eq:HS}, are indeed equivalent.

\medskip

\noindent\emph{Asymptotic behavior. } We use the  notation $p_R(|x|/R,t) := p_\infty(x,t)$. Now
$p_R$ is solved on the unit ball as
\begin{equation}
\label{eq:ode.tumor.spheroid}
 -p_R''(r) - \frac{N-1}r p_R'(r) = R^2
\Phi(p_R(r)) \quad \mbox{ for } \ 0<r<1,  \qquad p_R(1)=0, \quad
p_R'(0)=0.
\end{equation}
By the comparison principle, $p_R$ is increasing in $R$ and $p_R
\leq p_M$.  Passing to the limit in the equation we find that
$p_R\to p_M$ as $R    \to  \infty$.  It is easy to see that $p_R'(1)
=O(R)$,  because we can build  sub- and supersolutions of the form
$\bar p(r)= p_M\left(1- \frac{\textrm{e}^{\lambda
r}}{\textrm{e}^{\lambda}}\right)$ with $\lambda =CR$ for appropriate
values of $C$. By elliptic regularity we conclude that
$\|p_R'\|_{L^\infty(\mathbb{R}\times(0,T))}$ is of order $O(R)$ and,
since $p_R' \leq 0$, it is bounded in $L^1(\mathbb{R}\times(0,T))$.
As a consequence,  $ \|p_R'\|_{L^2(\mathbb{R}\times(0,T))}$ is of
order $O(\sqrt R)$.

Next we multiply equation \eqref{eq:ode.tumor.spheroid} by $p_R'$
and integrate in $(0,1)$ to obtain
$$
\frac {1}{2} \left(p_R'(1)\right)^2  + (N-1)
\int_0^1\frac{\left(p_R'(r)\right)^2}{r}\, dr = R^2 Q(p_R(0))
\approx R^2 Q(p_M).
$$
From the above $L^2$ estimate on $p_R' $, we conclude that $p_R'(1)
\approx R \sqrt{2 Q(p_M)}$. Going back to the dimensionalized
variable, \eqref{apts_speed} follows.

\medskip

\noindent\emph{One-dimensional profile. } We consider moving
coordinates,  $s=r-R(t)$, $\hat p_t(s)=p_\infty(s+R(t),t)$. Then
$\hat p_t(s)$ is positive if and only if $s<0$. Hence,
$$
\hat p_t''(s)+\frac{N-1}{s-R(t)}p_t'(s)+\Phi(p_t(s))=0\quad\text{for
} s<0, \qquad \hat p_t'(0)=-R'(t).
$$
Since $R(t)\to\infty$ as $t\to\infty$,  we can pass to the limit
uniformly on compact sets $-K\le s\le 0$ to obtain
$$
\hat p_\infty''(s)+\Phi(p_\infty(s))\quad \text{for } s<0, \qquad
\hat p_\infty'(0)=-\sqrt{2 Q(p_M)}.
$$

\qed

\noindent\emph{Remark. } The asymptotic growth
rate~\eqref{apts_speed} is generic for solutions to~\eqref{eq:HS},
since one can always  sandwich, after a certain delay, any \lq\lq
standard'' initial data between two particular profiles of the
form~\eqref{apts_rho}  which asymptotically travel with the same
speed~\eqref{apts_speed}.

\subsection{Examples with strong time discontinuities for the pressure}
\label{sec:timedisc}
Here we produce examples of solutions with a jump discontinuity in
the pressure as a function of time.

The first and simplest example is constructed in radial symmetry. It
consists of an initial datum which is constant $c_0<1$ in a ball
centered at $0$ and of radius 1  and zero otherwise. Then the
solution will have 0 pressure for a certain time, $0<t<t_1$, in
which the density grows exponentially according to
equation~\eqref{eq:HS}, so that
$$
\varrho_\infty(x,t)=c_0\textrm{e}^{\Phi(0)t}\mathds{1}_{B_1(0)}.
$$
At $t=t_1=-\frac{\log c_0}{\Phi(0)}$, $\varrho_\infty$ reaches the
threshold level $\varrho_\infty=1$ in $B_1$. It cannot happen that
$p_\infty(t)=0$ for $t>t_1$, since otherwise $\varrho_\infty$ would
violate the bound $\varrho_\infty \leq 1$. Therefore, $p_\infty(t)$
satisfies $\Delta p_\infty +\Phi(p_\infty)=0$ in some ball
$B_{R(t)}$ with $R(t)>0$. An easy argument using the maximum
principle shows that $R(t)\ge1$. It means that $p_\infty(t)$ is
equal or larger than the solution $\bar p_\infty$ of the problem
$\Delta \bar p_\infty +\Phi(\bar p_\infty)=0$, $\bar p_\infty=0$ for
$|x|=1$,  $\bar p_\infty>0$ for $|x|<1$. Taking the limit
$t\downarrow t_1$,  we obtain a jump discontinuity in time at
$t=t_1$ for every $|x|<1$. Note that $p_\infty(t)$ is a
discontinuous function of time with values in any
$L^p(\mathbb{R}^N)$, $p\ge1$.

The second example is best presented in one space dimension and
consists of the evolution of an initial data consisting of two
copies of the standard example presented in section
\ref{sec:spheroids} (tumor spheroid), which have now an interval as
support.  We locate the  supports at a distance from each other of
say 1. For a time interval, $0<t<t_1$ they evolve independently. At
$t=t_1$, the supports meet and an easy application of the maximum
principle, shows that the solution is strictly positive in the union
of the two intervals, therefore, larger than the solution
corresponding to the union of the two intervals at $t=t_1$. An easy
inspection of the solutions at times $t_1^-$ and $t_1^+$ shows that
there is a jump discontinuity in the space pressure profiles. This
example can be adapted to several space dimensions by replacing
disjoint intervals by disjoint concentric annuli.

\subsection{The effect of the equation on $\varrho_\infty$}
\label{subsect:effect.of.eq.density}

We consider the cell density  (defined
for  short enough times) with two discontinuities
$$
\varrho_\infty(t) =(1- q(t)) \mathds{1}_{B_{R_1(t)}} + q(t)
\mathds{1}_{B_{R_2}}, \qquad  R_1(t)< R_2,\quad q(t)\le 1.
$$
We claim that the correct dynamics is defined by the speed and
values of $q(t)$ given through
$$
R_1'(t)=\frac1{1- q(t)}  |\nabla p_\infty(R_1(t), t)|, \qquad  q'(t) = q(t) \Phi(0) \quad\text{while }q(t)\le1,
$$
with $p_\infty$ vanishing outside the ball of radius $R_1(t)$ and
satisfying $-\Delta p_\infty =\Phi(p_\infty)$ in $B_{R_1(t)}$, with zero Dirichlet boundary condition.

These rules are simply derived from equation \eqref{eq:HS}.  For
$R_1(t) < |x|< R_2$ we have $p_\infty=0$ and thus the equation
\eqref{eq:HS} is reduced to $\partial_t\varrho_\infty(t) = \Phi(0)\varrho_\infty(t)$  which
gives us the evolution of $q(t)$. For $|x|\le R_1(t)$, we have $ \partial_t \varrho_\infty = R_1'(t) (1- q(t))
\delta_{\{|x|=R_1(t)\} }$, while $\Delta p_\infty+\varrho_\infty\Phi(p_\infty) =|\nabla p_\infty|\delta_{\{|x|=R_1(t)\} }$. Using the equation
\eqref{eq:HS}, we get the dynamics for $R_1(t)$.

\

\noindent {\large\bf Acknowledgments.} Authors FQ and JLV  partially
supported by  Spanish project
MTM2011-24696. The paper was started while they where visiting UPMC.
They are indebted to this institution for the warm hospitality. BP
is supported by Institut Universitaire de France.

\vskip 1cm



\end{document}